\documentclass{article}

\usepackage[dvipsnames]{xcolor}
\usepackage{amsmath, amssymb, amsthm} 
\usepackage{geometry} 
\usepackage{graphicx} 
\usepackage{float} 
\usepackage{dsfont} 
\newcommand{\one}{\mathds{1}}
\usepackage{stmaryrd} 
\SetSymbolFont{stmry}{bold}{U}{stmry}{m}{n}   
\usepackage{nicefrac} 
\usepackage{comment}
\usepackage{upgreek}
\usepackage{relsize}

\usepackage{hyperref}

\usepackage{hyperref}
\hypersetup{
	linktocpage = true,
	colorlinks = true,
	linkcolor = {RoyalBlue},
	urlcolor = {RoyalBlue},
	citecolor = {Green}}

\usepackage{booktabs}
\usepackage{subcaption}
\DeclareFontShape{OT1}{cmr}{m}{scit}{<->ssub*cmr/m/sc}{} 

 \usepackage{placeins} 
 \usepackage{multirow}
\usepackage{pgfplots}
\usepackage{pgfplotstable}
\pgfplotsset{compat=1.18}   
\definecolor{seriesA}{HTML}{0072B2} 
\definecolor{seriesB}{HTML}{FF7E00} 
\definecolor{seriesC}{HTML}{8A2BE2} 
\newcommand{\xscalemode}{log} 

\usepackage{algorithm}      
\usepackage{algpseudocode}  

\newcommand{\resdir}{star_p150_gaussian}
\newcommand{\resdirB}{star_p150_bary}
 \newcommand{\mode}{small_majority}

\newcommand{\labU }{\texttt{`uniform'}{}}
\newcommand{\labLM}{\texttt{`large minority'}{}}
\newcommand{\labSM}{\texttt{`small majority'}{}}

\newcommand{\tabMedianeA}[2]{
        \begin{tabular}{@{}>{\centering\arraybackslash}m{1.2em}
                       *{6}{>{\centering\arraybackslash}m{0.13\textwidth}}@{}}\small
        & $\WBar_2$ & $\LWBar_2$ & $\WMed_1$ &  $\WMed_2$ (MM) & $\LWMed_2$ & $\QrLWm$ \\[2pt]
        \rotatebox{90}{healthy}
        & \fbox{\includegraphics[width=\linewidth]{images/#1/#2/WBar2_healthy.png}}
        & \fbox{\includegraphics[width=\linewidth]{images/#1/#2/LWBar2_healthy.png}}
        & \fbox{\includegraphics[width=\linewidth]{images/#1/#2/WMed1_Sink_healthy.png}}
        & \fbox{\includegraphics[width=\linewidth]{images/#1/#2/Wm2mm_healthy.png}}
        & \fbox{\includegraphics[width=\linewidth]{images/#1/#2/LWm2_healthy.png}}
        & \fbox{\includegraphics[width=\linewidth]{images/#1/#2/QrLWm2_healthy.png}} \\[2pt]
        \rotatebox{90}{corrupted}
        & \fbox{\includegraphics[width=\linewidth]{images/#1/#2/WBar2_corrupted.png}}
        & \fbox{\includegraphics[width=\linewidth]{images/#1/#2/LWBar2_corrupted.png}}
        & \fbox{\includegraphics[width=\linewidth]{images/#1/#2/WMed1_Sink_corrupted.png}}
        & \fbox{\includegraphics[width=\linewidth]{images/#1/#2/Wm2mm_corrupted.png}}
        & \fbox{\includegraphics[width=\linewidth]{images/#1/#2/LWm2_corrupted.png}}
        & \fbox{\includegraphics[width=\linewidth]{images/#1/#2/QrLWm2_corrupted.png}}
        \end{tabular}
}

\newcommand{\tabMedianeB}[2]{
        \begin{tabular}{@{}>{\centering\arraybackslash}m{1.2em}
                       *{7}{>{\centering\arraybackslash}m{0.11\textwidth}}@{}}\small
        & $\WBar_2$ & $\LWBar_2$ & $\WMed_1$ &  $\WMed_2$ (MM) &  $\WMed_2$ (irls)& $\LWMed_2$ & $\QrLWm$ \\[2pt]
        \rotatebox{90}{healthy}
        & \fbox{\includegraphics[width=\linewidth]{images/#1/#2/WBar2_healthy.png}}
        & \fbox{\includegraphics[width=\linewidth]{images/#1/#2/LWBar2_healthy.png}}
        & \fbox{\includegraphics[width=\linewidth]{images/#1/#2/WMed1_Sink_healthy.png}}
        & \fbox{\includegraphics[width=\linewidth]{images/#1/#2/Wm2mm_healthy.png}}
        & \fbox{\includegraphics[width=\linewidth]{images/#1/#2/Wm2irls_healthy.png}}
        & \fbox{\includegraphics[width=\linewidth]{images/#1/#2/LWm2_healthy.png}}
        & \fbox{\includegraphics[width=\linewidth]{images/#1/#2/QrLWm2_healthy.png}} \\[2pt]
        \rotatebox{90}{corrupted}
        & \fbox{\includegraphics[width=\linewidth]{images/#1/#2/WBar2_corrupted.png}}
        & \fbox{\includegraphics[width=\linewidth]{images/#1/#2/LWBar2_corrupted.png}}
        & \fbox{\includegraphics[width=\linewidth]{images/#1/#2/WMed1_Sink_corrupted.png}}
        & \fbox{\includegraphics[width=\linewidth]{images/#1/#2/Wm2mm_corrupted.png}}
        & \fbox{\includegraphics[width=\linewidth]{images/#1/#2/Wm2irls_corrupted.png}}
        & \fbox{\includegraphics[width=\linewidth]{images/#1/#2/LWm2_corrupted.png}}
        & \fbox{\includegraphics[width=\linewidth]{images/#1/#2/QrLWm2_corrupted.png}}
        \end{tabular}
}

\title{Quadratically Regularized Medians in Metric and \\ Wasserstein Spaces}
\author{Benoît Bonnet-Weill\footnote{Laboratoire des Signaux et Systèmes, Université Paris-Saclay, CNRS, CentraleSupélec, 91190 Gif-sur-Yvette, France. \textit{Email:} \texttt{benoit.bonnet-weill@centralesupelec.fr}}\;, Matteo Tacchi-Bénard\footnote{CNRS, GIPSA-lab, Univ. Grenoble Alpes, Grenoble, 38000, France}\; and Nicola Zaupa\footnote{CNRS, GIPSA-lab, Univ. Grenoble Alpes, Grenoble, 38000, France}}

\newcommand{\MB}[1]{\textcolor{red}{Matt: #1}}
\newcommand{\NZ}[1]{{\color{ForestGreen}{Nicola: #1}}}
\newcommand{\red}[1]{{\color{red}#1}}

\numberwithin{equation}{section}
\newtheorem{thm}{Theorem}[section]
\newtheorem{prop}[thm]{Proposition}
\newtheorem{cor}[thm]{Corollary}
\newtheorem{lem}[thm]{Lemma}
\theoremstyle{definition}
\newtheorem{dfn}[thm]{Definition}
\newtheorem*{dfnbis}{Definition}
\newtheorem{rem}[thm]{Remark}
\newtheorem{ex}[thm]{Example}

\newcommand{\scF}{\textsc{F}}
\newcommand{\Mean}{\textsc{Mean}}
\newcommand{\GMed}{\textsc{GMed}}
\newcommand{\QrGm}{\textsc{QrGm}}
\newcommand{\WBar}{\textsc{WBar}}
\newcommand{\WMed}{\textsc{WMed}}
\newcommand{\LWBar}{\textsc{LWBar}}
\newcommand{\LWMed}{\textsc{LWMed}}

\newcommand{\QrWm}{\textsc{QrWm}}
\newcommand{\QrLWm}{\textsc{QrLWm}}

\newcommand{\Eot}{\textsc{Eot}}
\newcommand{\SdBar}{\textsc{SdBar}}
\newcommand{\EBar}{\textsc{EBar}}

\newcommand{\CmD}{\textsc{CmD}}
\newcommand{\msupp}[2]{\mathbb{I}_{#1 \neq #2}}

\newcommand{\loc}{\mathrm{loc}}
\newcommand{\dist}{\mathrm{dist}}
\newcommand{\supp}{\mathrm{supp}}
\newcommand{\Leb}{\mathrm{Leb}}
\newcommand{\Id}{\mathrm{Id}}
\newcommand{\conv}{\mathrm{conv}}
\DeclareMathOperator*{\argmin}{argmin}
\DeclareMathOperator*{\argmax}{argmax}
\newcommand{\inte}[1]{\mathrm{int}(#1)}
\newcommand{\diam}[1]{\mathrm{diam}(#1)}
\newcommand{\spn}{\mathrm{span}}

\newcommand{\rmx}{\mathrm{x}}

\newcommand{\E}{\mathbb{E}}
\newcommand{\G}{\mathbb{G}}
\newcommand{\I}{\mathbb{I}}
\newcommand{\J}{\mathbb{J}}
\newcommand{\N}{\mathbb{N}}
\newcommand{\Pb}{\mathbb{P}}
\newcommand{\R}{\mathbb{R}}
\newcommand{\X}{\mathbb{X}}

\newcommand{\cC}{\mathcal{C}}

\newcommand{\cH}{\mathcal{H}}
\newcommand{\cG}{\mathcal{G}}
\newcommand{\cJ}{\mathcal{J}}
\newcommand{\cL}{\mathcal{L}}
\newcommand{\cM}{\mathcal{M}}
\newcommand{\cN}{\mathcal{N}}
\newcommand{\cP}{\mathcal{P}}
\newcommand{\cS}{\mathcal{S}}
\newcommand{\cX}{\mathcal{X}}
\newcommand{\cY}{\mathcal{Y}}
\newcommand{\cK}{\mathcal{K}}

\newcommand{\bT}{\boldsymbol{T}}
\newcommand{\bV}{\boldsymbol{V}}

\newcommand{\bfx}{\mathbf{x}}

\newcommand{\ba}{\boldsymbol{a}}
\newcommand{\bof}{\boldsymbol{f}}

\newcommand{\bp}{\boldsymbol{p}}
\newcommand{\bu}{\boldsymbol{u}}
\newcommand{\bx}{\boldsymbol{x}}
\newcommand{\by}{\boldsymbol{y}}
\newcommand{\bz}{\boldsymbol{z}}

\newcommand{\bdelta}{\boldsymbol{\delta}}
\newcommand{\blambda}{\boldsymbol{\lambda}}
\newcommand{\bnu}{\boldsymbol{\nu}}
\newcommand{\bsigma}{\boldsymbol{\sigma}}

\newcommand{\btheta}{\boldsymbol{\theta}}
\newcommand{\bgamma}{\boldsymbol{\gamma}}
\newcommand{\bmu}{\boldsymbol{\mu}}
\newcommand{\bvtheta}{\boldsymbol{\vartheta}}
\newcommand{\bG}{\boldsymbol{G}}
\newcommand{\bbeta}{\boldsymbol{\upbeta}}

\usepackage{dsfont}
\newcommand{\dsf}{\textnormal{\textsf{d}}}

\renewcommand{\rho}{\varrho}
\begin{document}

\maketitle

\begin{abstract}
In this paper, we introduce a new class of robust aggregators in Wasserstein spaces. After proving that all the commonly used Wasserstein barycenters are manipulable by one single corrupted input through a novel first-order optimality condition, we generalize both to the metric setting and to fractional exponents the concept of Lipschitz resilience introduced recently by Allouah, Guerraoui, Hoang and Villemaud. We proceed to show that adding a quadratic regularization to Wasserstein medians defined either in terms of the $W_1$-metric or linearized $W_2$-metric are resilient in that sense, by leveraging fine convexity estimates along generalized geodesics and flat interpolations. We finally show that iteratively solving the latter problem amounts to performing a Majorization-Minimization approximation of $W_2$-medians. 
\end{abstract}

\tableofcontents


\section{Introduction} \label{section:intro}

Whether to reach a collective decision, merge independent estimates or train a shared model, one often needs to summarize several pieces of information into a single representative one. In an increasing number of applications these pieces are no longer points but probability distributions, which encode the uncertainty, heterogeneity or noise of the sources they stem from. A natural way to aggregate them is the Wasserstein barycenter, introduced by Agueh and Carlier~\cite{Agueh2011} and made computationally tractable by Cuturi and Doucet~\cite{Cuturi2014}. It was quickly followed by successful applications such as texture mixing~\cite{Rabin2012} and shape interpolation~\cite{Solomon2015}. An overview of recent advances can be found in~\cite{PeyreComputational}, and, closer to our concerns, in~\cite{Panaretos2020} for statistics in Wasserstein spaces. In particular, the information to be aggregated might itself consist of posterior distributions of some parameters, computed in parallel on subsets of a large dataset:~\cite{Scott2016} averages Monte Carlo draws across data subsets, while~\cite{Srivastava2018} computes the Wasserstein barycenter of those posterior laws. In contrast,~\cite{Backhoff2022} replaces classical Bayesian predictive estimators by the barycenter of models weighted by their posterior law. However, like the arithmetic mean in Euclidean spaces, the barycenter of probability laws is not robust to outliers and adversarial attacks, and a single corrupted entry can shift it arbitrarily far, even in the presence of entropic regularizers (see Section~\ref{section:manipulability}).

A standard tool for assessing the vulnerability of data aggregators to outliers and attacks is the concept of \emph{breakdown point}~\cite{Hampel1971,Donoho1983,Huber2009}, which is defined as the minimal proportion of the input data one needs to control in order to be able to arbitrarily manipulate the output representation. For instance, the usual arithmetic mean has breakdown point $\nicefrac{1}{n}$ ($n$ being the number of input data), whereas the geometric median or Torricelli point, which minimizes the sum of distances to the inputs~\cite{Weiszfeld1937,Kemperman1987,Minsker2015}, achieves an asymptotic breakdown point (when $n$ goes to infinity) of $\nicefrac{1}{2}$, which is optimal among translation equivariant aggregators~\cite{Lopuhaa1991}. Medians have been extended from the initial Euclidean case to more general settings such as Riemannian manifolds~\cite{Fletcher2009}, but also spaces of probability measures, first to robustly combine subset posteriors~\cite{Minsker2017}, then as genuine Wasserstein medians~\cite{Carlier2024,You2025}. Other robust aggregators rely on discarding or down-weighing suspicious inputs. Popular examples include the median of means~\cite{Nemirovsky1983,Jerrum1986,Alon2002}, the trimmed mean~\cite{Tukey1963,Stigler1973}, and the Krum and multi-Krum rules~\cite{Blanchard2017,Bareilles2026}, which are widely used in Byzantine and federated learning~\cite{Yin2018,Pillutla2022,Guerraoui2024}, see also~\cite{Lugosi2019} for a survey in the Euclidean setting. We finally mention a complementary line of work which operates by robustifying the optimal transport distance itself rather than the aggregation rule~\cite{Mukherjee2021,Nietert2022,Le2021}. \medskip

While it provides fair robustness guarantees, the breakdown point criterion itself has important limitations, because it only assesses the possibility of \emph{arbitrary} manipulation in a binary fashion (manipulable or not manipulable), but never rules out strategic voting~\cite{ElMhamdi2023} and remains agnostic with respect to bounded manipulation. For instance, a minority of attackers can still significantly perturb the geometric median. While the existence of a closed form for the arithmetic mean contributes to it being easily inverted and hence manipulated, in contrast the geometric median is significantly more computationally expensive, as it requires to solve an optimization problem that can often be ill-conditioned (when the optimum is one of the input data points) or have a non-unique minimizer. That is why in this work, we study and extends a notion of quantitative robustness named Lipschitz resilience (together with an original Hölder variant), recently introduced in~\cite{Allouah2024}. The latter requires that the influence of any subset of inputs (be it a minority or a majority) on the aggregate be bounded by a power of their cumulated weights. This has two main advantages: first, resilient aggregators can never be arbitrarily manipulated, and second the manipulation bound grows smoothly with the weight of the attackers. Such strong robustness properties come at the price of loosing exactness and equivariance, which are two desirable representativeness properties of usual aggregators (see Section~\ref{section:Resilient}). Still, Lipschitz resilience is used in practice for robust collaborative filtering and content recommendation on the Tournesol platform~\cite{Hoang2021,Fageot2024}, and echoes the bounded sensitivity under differential privacy studied in~\cite{Dwork2006, Dwork2014}. 

In what follows, we introduce the prototypical example of Lipschitz resilient aggregator, whose properties will play a central role in Section \ref{subsection:QrLWm}. Therein, we assume that $(\cX,\|\cdot\|_{\cX})$ is a Banach space and let $(\cH, \| \cdot\|_\cH)$ be a Hilbert space such that $\cH \subset \cX$ and $\|\cdot\|_{\cX} \leq \|\cdot\|_\cH$.

\begin{dfnbis}[Quadratically regularized geometric median] 
Given a regularization parameter $L>0$ and a collection of weights $\blambda \in (\R_+^{\star})^n$, the \emph{quadratically regularized geometric median} with anchor $x_0 \in \cH$ of a family of points $\bx = (x_1,\ldots,x_n) \in \cH^n$ is defined by 
\begin{equation} 
\QrGm^{\blambda}_L(\bx \,|\,x_0) := \argmin_{y \in \cH} \frac{1}{L} \cdot \|y - x_0\|^2_\cH + \sum_{i=1}^n \lambda_i \cdot \|y - x_i\|_{\cX}.
\end{equation}
\end{dfnbis}

\noindent A small regularization parameter $L > 0$ pulls the geometric median towards the anchor $x_0$. From a bayesian viewpoint, it is the maximum a posteriori estimator for a Laplace likelihood and a Gaussian prior centered at $x_0$, with $L$ measuring the confidence placed in the prior (see Remark~\ref{rem:bayes}). From a variational viewpoint, it is the proximal operator of the median objective evaluated at $x_0$. Likewise, its Wasserstein counterpart is one step of the Jordan-Kinderlehrer-Otto scheme~\cite{Jordan1998,AGS} for the $W_1$-median functional studied in \cite{Carlier2024}. The strong convexity brought by the quadratic term is the key ingredient to both the resilience of these aggregators, and their favorable numerical behavior when computed using a regularized Weiszfeld algorithm.

Our analysis builds on several lines of work. Regularized Wasserstein barycenters have been extensively studied~\cite{Bigot2019, Cuturi2016, Janati2020, Carlier2021, Chizat2025, Genevay2020}, mostly from the prism of entropic regularization and Sinkhorn divergence~\cite{Feydy2019} and always with a mean-type data term, rather than median. Linearized optimal transport~\cite{Wang2013,Merigot2020,Moosmuller2023,Werenski2024}, together with the quantitative stability of transport maps~\cite{Delalande2023} and of barycenters~\cite{CarlierD2024}, turns Wasserstein aggregation into a Hilbertian problem, which we exploit for linearized medians. Finally, the fixed-point approach to Wasserstein barycenters of~\cite{Alvarez2016} can be recast as a majorization-minimization (MM) scheme~\cite{Hunter2004,Lange,Bolte2016}, which inspires an approximation scheme for un-regularized $W_2$-medians.

\paragraph{Contribution.} In a more sophisticated version of the above Hilbertian example, data could be corrupted by some noise with known (or estimated) distribution. They could also be subject to epistemic uncertainty if they are inferred from incomplete information, see e.g.~\cite{Allouah2024} on sparse voting. Moreover, the prior and posterior distributions may also be considered uncertain, so that the full aggregation problem would be meaningfully posed in a measure space. We propose to leverage Wasserstein metrics in order to directly aggregate probability measures (possibly built through non-parametric estimation or any meaningful process depending on the application) that may not agree with each other. In such a setting, Lipschitz and H\"older resilience would be instrumental in preventing adversarially generated measures from contaminating the aggregation too much, and finding a truly meaningful consensus. In this context, this article makes the following main contributions. 
\begin{enumerate}
    \item[$\diamond$] We extend the optimality of the $\nicefrac{1}{2}$ breakdown point~\cite{Lopuhaa1991,Carlier2024} to aggregators on metric spaces that are equivariant under an unbounded group of isometries (Proposition~\ref{prop:breakdown}).
    \item[$\diamond$] We show that Wasserstein barycenters, as well as their (doubly) entropic, Sinkhorn and linearized variants, have the same breakdown point as the Euclidean arithmetic mean, although they are not exactly manipulable in general (Section~\ref{section:manipulability}).
    \item[$\diamond$] We provide a general criterion for Hölder and Lipschitz resilience based on uniform convexity estimates, and show that $\QrGm$ is Lipschitz resilient (Section~\ref{section:Resilient}).
    \item[$\diamond$] We introduce two quadratically regularized medians in Wasserstein spaces. The first one, called $\QrWm$, is built on $W_1$-medians and shown to satisfy handy moment support bounds, and to be H\"older resilient (Section~\ref{subsection:QrWm}). The second one, coined $\QrLWm$, involves linearized $W_2$-distances. It is shown to be Lipschitz resilient, and to admits a simple Euclidean representation together with density bounds (Section~\ref{subsection:QrLWm}).
    \item[$\diamond$] We propose a majorization-minimization scheme to approximate $W_2$-medians based on iterating $\QrLWm$ from an initial prior $\mu_0$, and show that its cluster points are critical for the vanilla median problem (Section~\ref{subsection:MajorizeMinimize}).
    %
\end{enumerate}

\section{Preliminaries and elementary results}
\label{section:Preli}

In what follows, we introduce preliminary material and notations about optimal transport and abstract aggregators over metric spaces. 


\subsection{Integration and optimal transport.}

In this section, we recall basic notions of integration and optimal transports over $\R^d$. Throughout, we fix some real number $p \in [1,+\infty)$ and let $(\cX,\dist)$ be a complete separable metric space. 

We shall denote by $\cM_{\loc}^+(\cX)$ and $\cM^+(\cX)$ the sets of positive and positive finite Radon measures on $\cX$, and write $\cP(\cX)$ for the subspace of Radon probability measures. Given any $\mu \in \cM^+(\cX)$, we let $L^p(\cX,\cY;\mu)$ stand for the space of $p$-Lebesgue integrable map valued in some finite dimensional Banach space $(\cY,\|\cdot\|_{\cY})$. We will say that a measure is absolutely continuous, which we write as $\mu \in \cP^{\mathrm{ac}}(\R^d)$, provided that $\mu \ll \Leb$. We lastly denote by $\cP_c^{\infty}(\R^d) \subset \cP^{\mathrm{ac}}(\R^d)$ the set of absolutely continuous probability measures with convex compact support, and whose density is both lower- and upper-bounded. In what follows, we define the set of \textit{transport plans} between two measures $\mu,\nu \in \cP(\cX)$ by 
\begin{equation*}
\Gamma(\mu,\nu) := \bigg\{ \gamma \in \cP(\cX^2) \; \Big| \; \pi^1_{\sharp} \gamma = \mu ~\text{and}~ \pi^2_{\sharp} \gamma = \nu \bigg\}, 
\end{equation*}
where ``$\sharp$'' is the standard pushforward operation, and $\pi^1,\pi^2 : \cX^2 \to \cX$ are the projections onto the first and second factor respectively. In what follows, we let 
\begin{equation}
\cP_p(\cX) := \left\{ \mu \in \cP(\cX) \; \Big| \; \int_{\cX} \dist(x_0,x)^p \, \mathrm{d}\mu(x)  < +\infty \right\}
\end{equation}
for some (and thus all) $x_0 \in \cX$, and denote by $\cP_c(\cX)$ the subset of probability measures with compact support. We also recall the definition of the \textit{Kullback-Leibler divergence}  
\begin{equation}
\mathrm{KL}(\mu \,|\, \nu) = \left\{
\begin{aligned}
& \int_{\cX^2} \log\left(\frac{\mathrm{d}\mu}{\mathrm{d}\nu}(x) \right) \, \mathrm{d}\mu(x) ~~ & \text{if $\mu \ll \nu$}, \\
& + \infty ~~ & \text{otherwise}, 
\end{aligned}
\right.
\end{equation}
for every $\mu,\nu \in \cM^+_{\loc}(\cX)$. Finally, in the particular case in which $(\cX,\|\cdot\|_{\cX})$ is a finite-dimensional Banach space, we denote by 
\begin{equation*}
\cM_p(\mu) := \bigg( \int_{\cX} \|x\|^p_{\cX} \, \mathrm{d} \mu(x) \bigg)^{\nicefrac{1}{p}}   
\end{equation*}
the moment of order $p \in [1,+\infty)$ of an element $\mu \in \cP_p(\cX)$. 

\begin{dfn}[Wasserstein distances and spaces]
The \emph{$p$-Wasserstein distance} between $\mu,\nu \in \cP_p(\R^d)$ is the quantity defined by 
\begin{equation} \label{eq:wasserstein_distance}
W_p(\mu,\nu) = \min \bigg\{ \bigg( \int_{\cX^2} \dist(x,y)^p \, \mathrm{d}\gamma(x,y) \bigg)^{\nicefrac{1}{p}}\; \Big| \; 
\gamma \in \Gamma(\mu,\nu) \bigg\}
\end{equation}
The set of minimizers of~\eqref{eq:wasserstein_distance} -- which are called \emph{optimal transport plans} -- is denoted by $\Gamma_o(\mu,\nu)$, and the $p$-Wasserstein space then refers to the metric space $(\cP_p(\cX),W_p)$. 
\end{dfn}

It is well-known that $(\cP_p(\cX),W_p)$ is a complete separable metric space, see e.g. \cite[Chapter 7]{AGS}. In the sequel, we recollect the definitions and some basic properties of several classical generalizations of the optimal transport problem $(\cX,\|\cdot\|_{\cX})$ is a Banach space.    

\begin{dfn}[Entropic optimal transport]
The \emph{entropic $p$-optimal transport problem} with regularization parameter $\alpha \in \R_+$ between two measures $\mu,\nu \in \cP_p(\cX)$ is defined by
\begin{equation}
\label{eq:EntropicOT}
\Eot_{p}(\mu,\nu \,|\, \alpha) := \min \bigg\{\int_{\cX^2} \dist(x,y)^p \, \mathrm{d}\gamma(x,y) + \alpha\cdot\mathrm{KL}(\gamma \,|\, \mu\otimes\nu) \; \Big| \; \gamma \in \Gamma(\mu,\nu) \bigg\} 
\end{equation}
The set of minimizers of~\eqref{eq:EntropicOT} -- which are called \emph{optimal entropic plans} -- is denoted by $\Gamma_o^{\alpha}(\mu,\nu)$
\end{dfn}

By definition of the entropic optimal transport problem, it clearly holds that 
\begin{equation*}
W_p^p(\mu,\nu) \leq \Eot_p(\mu,\nu\,|\,\alpha) 
\end{equation*}
for every $\mu,\nu \in \cP_p(\R^d)$, and likewise $\Eot_{p,0}(\mu,\nu) = W_p^p(\mu,\nu)$. Moreover, the quantity $\Eot(\cdot,\cdot\,|\,\alpha)$ satisfies a (pseudo) triangle inequality and slack upper bound in terms of the moments of the inputs. These are folklore results that we state in the following lemma for the sake of self-containedness. 

\begin{lem}[Pseudo-triangle inequality and crude upper bound for $\Eot_p$]
\label{lem:PseudoTriangle}
Given any $p \in [1,+\infty)$ and $\alpha \in \R_+$, it holds that
\begin{equation}
\label{eq:PseudoTriangle}
\Eot_p(\mu,\sigma \,|\,\alpha) \leq 2^{p-1} \cdot \Big( \Eot_p(\mu,\nu \,|\,\alpha) + \Eot_p(\nu,\sigma \,|\,\alpha) \Big)
\end{equation}
as well as
\begin{equation}
\label{eq:CrudeBound}
\Eot_p(\mu,\nu \,|\,\alpha) \leq 2^{p-1} \cdot \Big( \cM_p^p(\mu) +\cM_p^p(\nu) \Big)
\end{equation}
for every $\mu,\nu,\sigma \in \cP_p(\cX)$
\end{lem}

The last variant of classical optimal transport that we will use in the sequel is the so-called \textit{linearized Wasserstein distance} over $(\cX,\|\cdot\|_{\cX}) := (\R^d,\|\cdot\|)$, which was first introduced in \cite{Wang2013} and studied quite extensively in recent works investigating the stability of optimal transport maps (see e.g. \cite{CarlierD2024,Delalande2023} and subsequent works). 

\begin{dfn}[Linearized Wasserstein distance]
For every $p > 1$, the \emph{linearized $W_p$-metric} around $\rho \in \cP_p^{\mathrm{ac}}(\R^d)$ is defined by 
\begin{equation} 
\label{eq:LOT}
W_{p,\rho}(\mu,\nu) = \min \Bigg\{ \bigg( \int_{\R^{3d}} \|y-z\|^p \, \mathrm{d}\bgamma(x,y,z) \bigg)^{\nicefrac{1}{p}} \; \Big| \; (\pi^1,\pi^2)_{\sharp} \bgamma \in \Gamma_o(\rho,\mu) ~\text{and}~ (\pi^1,\pi^3)_{\sharp} \bgamma \in \Gamma_o(\rho,\nu) \Bigg\}
\end{equation}
for every $\mu,\nu \in \cP_p(\R^d)$. Moreover, it clearly holds that
\begin{equation}
\label{eq:LinearizedWassIneq}
W_p(\mu,\nu) \leq W_{p,\rho}(\mu,\nu)
\end{equation}
since $(\pi^2,\pi^3)_{\sharp} \bgamma \in \Gamma(\mu,\nu)$ for every plan $\bgamma \in \cP(\R^{3d})$ as above.
\end{dfn}

Linearized Wasserstein distances are most commonly considered in the case where $p=2$, as it follows then from Brenier's theorem (see e.g. \cite[Theorem 1.22]{OTAM}) that the latter reduce to Hilbertian distances between (uniquely determined) optimal transport maps, which are themselves in one-to-one correspondence with the convex cone $\nabla \cC(\R^d)$ of gradients of convex function. These handy facts are recalled in the following theorem.

\begin{thm}[Linearized Wasserstein distance around regular measures] 
\label{thm:Brenier}
Suppose that $p=2$ and let $\rho \in \cP_2^{\mathrm{ac}}(\R^d)$. Then $\Gamma_o(\rho,\mu) = \{(\Id,T_{\rho}^{\mu})_{\sharp} \rho\}$ for a uniquely determined $T_{\rho}^{\mu} \in L^2(\R^d,\R^d;\rho)$, and 
\begin{equation}
\label{eq:LinearizedEquality}
W_{2,\rho}(\mu,\nu) = \|T_{\rho}^{\nu} - T_{\rho}^{\mu} \|_{L^2(\rho)}.
\end{equation}
Moreover, a map $T_{\rho}^{\mu} \in L^2(\R^d,\R^d;\rho)$ is such that $(\Id,T_{\rho}^{\mu})_{\sharp} \rho \in \Gamma_o(\rho,\mu)$ \emph{if and only if} $T_{\rho}^{\mu} = \nabla \varphi$ for some convex function $\varphi : \R^d \to \R \cup \{+\infty\}$. 
\end{thm}

We end this preliminary section by recalling the expression of the first-order variation of the squared Wasserstein distance at an absolutely continuous measure, which is borrowed from \cite[Corollary 10.2.7]{AGS}. 

\begin{thm}[Fréchet differentiability of the squared $W_2$-metric]
\label{thm:WassGrad}
Suppose that $p=2$ and let $\nu \in \cP_2(\R^d)$ be arbitrary. Then for each $\mu \in \cP_2^{\mathrm{ac}}(\R^d)$, one has that 
\begin{equation*}
\tfrac{1}{2} W_2^2(\nu,\sigma) = \tfrac{1}{2} W_2^2(\nu,\mu) - \big\langle T_{\mu}^{\sigma} - \Id , T_{\mu}^{\nu} - \Id \big\rangle_{L^2(\mu)} + \, o(W_2(\mu,\sigma))
\end{equation*}
for every $\sigma \in \cP_2(\R^d)$
\end{thm}

Following the above result, one may define the \emph{Wasserstein gradient} of the squared Wasserstein distance from some $\nu \in \cP_2(\R^d)$ at $\mu \in \cP_2^{\mathrm{ac}}(\R^d)$ as the map $\nabla_{\mu} W_2^2(\nu,\mu) \in L^2(\R^d,\R^d;\mu)$ given by  
\begin{equation}
\label{eq:GradWass}
\nabla_{\mu} W_2^2(\nu,\mu) := 2(\Id - T_{\mu}^{\nu})
\end{equation}
for every $\mu \in \cP_2^{\mathrm{ac}}(\R^d)$. 


\subsection{Abstract aggregators and their breakdown point}

In this section, we introduce several basic notions about abstract aggregators in metric spaces. We also revisit some known results about Fréchet means, medians, and their breakdown points in this general setting.  

\begin{dfn}[Mismatch support]
Any set $\cS$ can be equipped with the discrete distance 
\begin{equation*}
\delta : (x,y) \mapsto
\left\{
\begin{aligned}
& 1 \quad & \text{if $x \neq y$}, \\
& 0 \quad & \text{otherwise.}
\end{aligned}
\right.
\end{equation*}
Consequently for $n \in \N^\star$ and $\blambda = (\lambda_1,\ldots,\lambda_n) \in \Delta^n$, we may equip $\cS^n$ with the weighted \emph{cumulative mismatch pseudo\footnote{If the $\lambda_i$ are all positive, then $\CmD_{\blambda}$ is a proper distance, but otherwise it may fail to be separated.} distance}
\begin{equation} 
\label{eq:cmd}
\CmD_{\blambda}  : (\bx,\by) \longmapsto \displaystyle\sum_{i=1}^n \lambda_i \cdot \delta(x_i,y_i).
\end{equation}
Then we define the \emph{mismatch support} of two elements $\bx,\by \in \cS^n$ by 
\begin{equation} \label{eq:mismatch}
\msupp{\bx}{\by} := \supp(\delta(\bx,\by)) = \big\{i \in \llbracket 1,n \rrbracket \,|\, x_i \neq y_i \big\},
\end{equation}
and note that $\CmD_{\blambda}$ is by construction the sum of the mismatch weights, namely 
\begin{equation*}
\CmD_{\blambda}(\bx,\by) = \sum_{i\in\msupp{\bx}{\by}} \lambda_i.
\end{equation*}
\end{dfn}

In the sequel, we let $(\cX,\dist)$ be a metric space and use the same notation to refer to the induced set distance $\dist(A,B) := \inf_{(a,b) \in A \times B} \dist(a,b)$. We also fix an integer $n \in \N^\star$ 

\begin{dfn}[Abstract weighted aggregator]
\label{def:Aggregators}
A set-valued mapping $\scF : (\R_+^{\star})^n \times \cX^n \rightrightarrows\cX$ is called a \emph{weighted aggregator} if it satisfies the following axioms.
\begin{enumerate}
\item[$(i)$] Data with null weights have no influence on $\scF^{\blambda}$, namely if $\lambda_i = 0$ and $\bx,\by \in \cX^n$ are such that $\msupp{\bx}{\by} \subset \{i\}$, then $\scF^{\blambda}(\bx) = \scF^{\blambda}(\by)$.
\item[$(ii)$] Data with nonzero weights have an influence on $\scF^{\blambda}$, namely for $i \in \llbracket 1,n \rrbracket$, there exists $\blambda \in \Delta^n$ with $\lambda_i > 0$ and $\bx,\by \in \cX^n$ with $\msupp{\bx}{\by} = \{i\}$ such that $\scF^{\blambda}(\bx) \neq \scF^{\blambda}(\by)$.
\item[$(iii)$] The mapping $\blambda \in (\R_+^{\star})^n \rightrightarrows \scF^{\blambda}(\bx) \subset $ is \emph{monotone} for every $\bx \in \cX^n$, in the sense that if given some $i \in \llbracket 1,n \rrbracket$ we take $\blambda,\blambda' \in (\R_+^{\star})^n$ such that $\I_{\blambda \neq \blambda'} \subset \{i\}$ and $\lambda_i \leq \lambda_i'$, then 
\begin{equation*}
\dist \Big( x_i , \scF^{\blambda'}(\bx) \Big) \leq \dist \Big( x_i , \scF^{\blambda}(\bx) \Big)
\end{equation*}
%
%
In words, increasing the weight of a fixed data point $x_i \in \cX$ decreases the distance of the aggregator to this data point.
%
\end{enumerate}
In the sequel, we will almost systematically consider aggregators whose weights $\blambda \in \Delta^n$ belong to the $n$-dimensional simplex. 
\end{dfn}

\begin{rem}[Defining aggregators on the positive orthant]
Our choice to provide an axiomatic for aggregators on the orthant $(\R_+^{\star})^n$ rather than on the simplex $\Delta^n$ directly is due to the fact that some properties, notably monotonicity, may be true on the former but fail on the latter. This is the case for the quadratically regularized medians we introduce below, which satisfy the aggregator axiomatic on $(\R_+^{\star})^n$ but are not monotone on $\Delta^n$.  
\end{rem}

\begin{ex}[Basic aggregators in metric spaces]
The \emph{Fr\'{e}chet median}, defined by 
\begin{subequations}
\begin{equation} \label{eq:fmed}
\textsc{FrMed}^{\blambda}(\bx) := \argmin_{y \in \cX} \sum_{i=1}^n \lambda_i \, \dist(x_i,y) \subset \cX
\end{equation}
is a weighted aggregator, so is the usual \emph{Fr\'{e}chet mean} 
\begin{equation} \label{eq:fmean}
\textsc{FrMean}^{\blambda}(\bx) := \argmin_{y \in\cX} \sum_{i=1}^n \lambda_i \, \dist(x_i,y)^2 \subset \cX.
\end{equation}
\end{subequations}
\end{ex}

\begin{dfn}[Breakdown point of an aggregator] 
Let $\blambda \in \Delta^n$ and assume that $(\cX,\dist)$ has infinite diameter, i.e. for every $R > 0$, there exist $x,y \in \cX$ such that $ \dist(x,y) \geq R$. Then, the \emph{breakdown point} of a weighted aggregator $\scF^{\blambda} : \cX^n \rightrightarrows \cX$ at an input $\bx \in \cX^n$ is defined by 
\begin{equation} \label{eq:breakdown}
\bbeta^n_{\dist}(\scF^{\blambda}(\bx)) := \inf \bigg\{\beta \in [0,1] \; \Big| \; \sup_{\CmD_{\blambda}(\bx,\by) \leq \beta} \dist_{H}\left(\vphantom{\sum} \scF^{\blambda}(\bx), \scF^{\blambda}(\by)\right) = +\infty \bigg\},
\end{equation}
where 
\begin{equation*}
\dist_H(A,B) = \max\bigg\{\sup_{a\in A} \dist(a,B), \; \sup_{b\in B} \dist(b,A) \bigg\} 
\end{equation*}
is the usual Hausdorff distance. Concretely, $\bbeta^n_{\dist}(\scF^{\blambda}(\bx))$ is the minimum proportion of mismatch required to arbitrarily shift $\scF^{\blambda}(\bx)$, measured in terms of the Hausdorff distance. When applicable, we will also define the \textit{asymptotic breakdown point}
\begin{equation*}
\bbeta^{\infty}_{\dist}(\scF(\bx)) = \lim\limits_{n\to\infty} \bbeta^n_{\dist}(\scF^{\bu_n}(\bx))
\end{equation*}
where $\bu_n = \nicefrac{\one_n}{n} := (\nicefrac{1}{n},\ldots,\nicefrac{1}{n}) \in \Delta^n$ is the uniform simplex vector.
\end{dfn}

\begin{ex}[Breakdown point of the arithmetic mean]
If $\cX$ is a vector space, then one can define the \emph{arithmetic mean} with weight $\blambda \in \Delta^n$ as 
\begin{equation} \label{eq:mean}
\Mean^{\blambda}(\bx) := \sum_{i=1}^n \lambda_i \, x_i.
\end{equation}
It is well-known that, if the distance comes from a norm, namely $\dist(\bx,\by) = \|\bx-\by\|$, then 
\begin{equation*}
\bbeta^n_{\|\cdot\|}(\Mean^{\blambda}(\bx)) = \min\big\{\lambda_i \,|\, i \in \llbracket 1,n \rrbracket, \; \lambda_i > 0\big\}. 
\end{equation*}
Indeed, choosing $j\in \argmin\{\lambda_i \,|\, i \in \llbracket 1,n \rrbracket, \; \lambda_i > 0\}$ and $\bx,\by \in \cX^n$ such that $y_i = x_i$ for $i \neq j$ yields
\begin{equation*}
\begin{aligned}
\dist\left(\vphantom{\sum} \Mean^{\blambda}(\bx) , \Mean^{\blambda}(\by) \right) & = \bigg\| \sum_{i=1}^n \lambda_i (x_i - y_i) \bigg\| \\
& = \lambda_{j} \|x_{j} - y_{j}\| \underset{\|y_{j}\|\to\infty}{\longrightarrow} +\infty. 
\end{aligned}
\end{equation*}
%
This allows to recover the commonly know facts that $\bbeta^n_{\|\cdot\|}(\Mean^{\bu_n}) = \nicefrac{1}{n}$ and $\bbeta^{\infty}_{\|\cdot\|}(\Mean) = 0$.
\end{ex}

\begin{ex}[Fréchet medians in geodesic spaces]
Suppose that $\cX$ is a geodesic space, namely for all $x_0,x_1 \in \cX$, there exists a curve $x(\cdot) \in \cC([0,1],\cX)$ with $x(0) = x_0$, $x(1)=x_1$ and a $C>0$ such that every $\tau \in [0,1]$ admits a neighborhood $I \subset [0,1]$ satisfying
\begin{equation*}
\dist(x(s),x(t)) = C|t-s|
\end{equation*}
whenever $s,t \in I$. Then, it was proven in~\cite{Carlier2024} that
\begin{equation*} 
\bbeta^n_{\dist}(\textsc{FrMed}^{\blambda}(\bx)) = \min \bigg\{ \sum_{j \in J} \lambda_j \; \Big| \; J \subset \llbracket 1,n \rrbracket ~~\text{and}~~ \sum_{j\in J} \lambda_j \geq \frac{1}{2} \bigg\}, 
\end{equation*}
and in particular, one recovers the famous results
\begin{equation*}
\bbeta^n_{\dist}(\textsc{FrMed}^{\bu_n}(\bx)) = \frac{\left\lfloor \nicefrac{(n+1)}{2} \right\rfloor}{n} \qquad \text{and} \qquad \bbeta^{\infty}_{\dist}(\textsc{FrMed}) = \frac{1}{2}. 
\end{equation*}
\end{ex}

It is well-known that in the Euclidean setting, the (asymptotic) breakdown point of the Fréchet median is optimal in the sense that it is the biggest reachable one within the class of translation-invariant aggregators, see~\cite{Lopuhaa1991}. This still holds in the metric case, provided there exists a group a group $(\G,\circ)$ of unbounded isometries whose elements act on $(\cX,\dist)$. Concretely, this means that every $g \in \G$ acts on $\cX$ through a group operation of the form $g \cdot x \in \cX$, with   
\begin{equation*}
\dist(g \cdot x,g\cdot y) = \dist(x,y)
\end{equation*}
for every $g \in \G$ and all $x,y \in \cX$. Requiring that $(\G,\circ)$ be unbounded means that for each $x \in \cX$ and $R>0$, there must exist $g^x_R \in \G$ such that
\begin{equation*}
\dist(x, g^x_R\cdot x) \geq R.
\end{equation*}
Such isometries play the role of translations in the metric reformulation of the proof in~\cite{Lopuhaa1991}, which is the object of the following proposition.

\begin{prop}[On the optimality of the $\nicefrac{1}{2}$-breakdown point] \label{prop:breakdown}
Suppose that $\scF^{\blambda} : \cX^n \rightrightarrows \cX$ is equivariant under the action of a group of unbounded isometries $(\G,\circ)$, namely
\begin{equation*}
\scF^{\blambda}(g\cdot\bx) = g\cdot \scF^{\blambda}(\bx)
\end{equation*}
for each $g \in \G$, where $g \cdot (x_1,\ldots,x_n) := (g\cdot x_1, \ldots, g\cdot x_n)$ and $g\cdot \scF^{\blambda}(\bx) = \{g\cdot \by \mid \by \in \scF^{\blambda}(\bx)\}$. Furthermore, assume that $\scF^{\blambda}(\bx)$ is nonempty and bounded for all $\bx \in \cX^n$. Then, it holds that 
\begin{equation}
\label{eq:RemBreakdown}
\sup\bigg\{\dist_H\Big(\scF^{\bu_n}(\bx), \scF^{\bu_n}(\by) \Big) \; \Big| \; \CmD_{\one_n}(\bx,\by) = |\msupp{\bx}{\by}| \leq \big\lceil \nicefrac{n}{2} \big\rceil \bigg\} = +\infty,
\end{equation}
which entails in particular that
\begin{equation*}
\bbeta^n_{\dist}(\scF^{\bu_n}(\bx)) \leq \frac{\lfloor \nicefrac{(n+1)}{2}\rfloor}{n} \underset{n\to\infty}{\longrightarrow} \frac{1}{2}
\end{equation*}
for every $\bx \in \cX^n$. 
\end{prop}

\begin{proof}
The proof goes by contradiction. Indeed, assume that the quantity in \eqref{eq:RemBreakdown} is finite and upper-bounded by some $M > 0$. Then, taking $\theta \in \scF^{\bu_n}(\bx)$ and for an arbitrary $\bx \in \cX^n$ setting $D := \diam{\scF^{\bu_n}(\bx)}$ and $R = 4M + D + 1$, the collections of points defined by 
\begin{equation*}
\by := \Big(x_1, \ldots, x_{\lfloor \nicefrac{n}{2} \rfloor}, g_R^{\theta}\cdot x_{\lfloor\nicefrac{n}{2}\rfloor+1}, \ldots, g_R^{\theta} \cdot x_n \Big) \in \cX^n  \qquad \text{and} \qquad \bz := (g_R^{\theta})^{-1} \cdot \by \in \cX^n
\end{equation*}
possess at most $\lfloor \nicefrac{(n+1)}{2} \rfloor$ mismatches with $\bx$, so that
\begin{equation*}
\dist_H \Big(\scF^{\bu_n}(\bx),\scF^{\bu_n}(\by) \Big) \leq M \qquad \text{and} \qquad \dist_H \Big(\scF^{\bu_n}(\bx),\scF^{\bu_n}(\bz) \Big) \leq M.
\end{equation*}
Moreover, it follows from the definition of Hausdorff distance that
\begin{equation*}
\dist(\theta, g_R^\theta \cdot \theta) \leq \dist\Big(\scF^{\bu_n}(\bx) , g_R^\theta\cdot\theta \Big) + D \leq \dist_H\Big(\scF^{\bu_n}(\bx), g_R^\theta\cdot \scF^{\bu_n}(\bx) \Big) + D.
\end{equation*}
Then by the triangle inequality, we know that $\dist_H(\scF^{\bu_n}(\by),\scF^{\bu_n}(\bz)) \leq 2M$, whence
\begin{align*}
4M+1 & = R - D  \\
& \leq \dist(\theta,g_R^\theta\cdot \theta) - D \\
& \leq \dist_H\Big( \scF^{\bu_n}(\bx), g_R^\theta \cdot \scF^{\bu_n}(\bx) \Big) \\
& \leq \underset{\leq M}{\underbrace{\dist_H\Big( \scF^{\bu_n}(\bx), \scF^{\bu_n}(\bz) \Big)}} + \underset{\leq 2M}{\underbrace{\dist_H\Big( \scF^{\bu_n}(\bz),\scF^{\bu_n}(\by) \Big)}} + \dist_H\Big( \scF^{\bu_n}(\by), g_R^\theta\cdot \scF^{\bu_n}(\bx) \Big) \\
& \leq 3 M + \dist_H\Big( \scF^{\bu_n}(\by), g_R^\theta\cdot \scF^{\bu_n}(\bx) \Big).
\end{align*}
This leads to a contradiction by noticing that $\scF^{\bu_n}(\by) = \scF^{\bu_n}(g_R^\theta \cdot \bz)$, since then by group action equivariance and isometry, it holds 
\begin{equation*}
\begin{aligned}
\dist_H\Big( \scF^{\bu_n}(\by), g_R^\theta\cdot \scF^{\bu_n}(\bx) \Big) & = \dist_H \Big( g_R^\theta\cdot\scF^{\bu_n}(\bz), g_R^\theta\cdot \scF^{\bu_n}(\bx) \Big) \\
& = \dist_H\Big( \scF^{\bu_n}(\bz), \scF^{\bu_n}(\bx) \Big) \\
& \leq M.   
\end{aligned}
\end{equation*} 
In conclusion, the quantity in \eqref{eq:RemBreakdown} cannot be finite. 
\end{proof}


\begin{dfn}[Geometric median]
The \emph{geometric median} of a collection of inputs $\bx \in \cX^n$ with weights $\lambda \in \Delta^n$ is defined by 
\begin{equation} \label{eq:GM}
\GMed^{\blambda}(\bx) := \argmin_{y \in \cX} \sum_{i=1}^n \lambda_i \, \|y - x_i\|.
\end{equation}
It is a special instance of the Fréchet median, and shares the same breakdown point. 
\end{dfn}

\begin{rem}[Concerning the robustness of geometric medians]
$\GMed^{\blambda}$ is obviously much more robust than $\Mean^{\blambda}$, since it can only be arbitrarily manipulated by a majority of attackers. Even so, a large minority of outliers is still not harmless, and can significantly perturb the aggregation. Supporting this statement, it was proved in \cite{ElMhamdi2023} that geometric medians in Euclidean spaces are not strategy-proof in general, in the sense that for each $k \in \N^\star$, there exist $\bx_k,\by_k \in (\R^d)^n$ such that $\msupp{\bx_k}{\by_k} = \{1\}$ and
\begin{equation*}
\|\GMed^{\blambda}(\by_k) - x_{k,1}\| \leq \frac{1}{k} \|\GMed^{\blambda}(\bx_k) - x_{k,1}\|.    
\end{equation*}
In other words, for any arbitrary multiplicative gain $k$, there exist configurations where lying on $x_1$ (i.e. replacing it with $y_1$) allows to divide by $k$ the distance to $x_1$. More generally, the fact that only a majority of attackers can \emph{arbitrarily} manipulate $\GMed^{\blambda}$ does not provide information on the impact a minority of attackers can have, except that it is bounded. 
\end{rem}

\begin{rem}[Concerning the computation cost of geometric medians]
Another problem with $\GMed^{\blambda}$ that we only touch upon in this contribution is its computational cost. Indeed, even in high dimensional spaces, computing $\Mean^{\blambda}$ is almost free, as it is simply a linear combination of the inputs with known coefficients (the weights). In contrast, computing $\GMed^{\blambda}$ requires to solve an optimization problem that is convex -- but not strongly --, which can be costly if the aggregation is done in real time. 
Still, besides robustness, another motivation for robust statistics is arguably to find aggregation techniques with robustness guarantees similar to $\GMed^{\blambda}$, while keeping the computational burden as low as possible.
\end{rem}


\section{Manipulable and robust aggregators in metric spaces}

In this section, we first prove that all commonly used barycenter-type aggregators in Wasserstein spaces are manipulable by a single adversarial input, and then develop a framework for resilient aggregators on metric spaces inspired by \cite{Allouah2024}.


\subsection{Manipulability of Wasserstein Barycenters} \label{section:manipulability}

In this section, we show that both classical, entropic, and Sinkhorn divergence barycenters in Wasserstein spaces have the same breakdown point measured in Wasserstein metrics as the usual Euclidean one, even with an entropic regularization. To the best of our knowledge, such a general result was not covered by the existing literature, and its proof mechanism relying on a first-order intrinsic optimality condition is likewise new.   

The seminal paper \cite{Agueh2011} put forth a notion of barycenter in Wasserstein spaces when $p \in (1,+\infty)$, which has since then become one of the choice methods for aggregating collections of probability measures. 

\begin{dfn}[Classical Wasserstein barycenters]
Given a collection of weights $\blambda \in \Delta^n$, the \emph{$W_p$-barycenter} of a family of measures $\bnu \in \cP_p(\R^d)^n$ is defined by 
\begin{equation}
\textsc{WBar}_p^{\blambda}(\bnu) = \argmin_{\mu \in \cP_p(\R^d)} \sum_{i=1}^n \lambda_i \, W_p^p(\nu_i,\mu). 
\end{equation}
\end{dfn}
Various generalizations of the usual Wasserstein barycenters involving entropy terms have been proposed in the literature, see e.g. \cite{Cazelles2018,Carlier2021}. Here, we focus on the \emph{doubly regularized entropic} framework investigated in \cite{Chizat2025}, as it subsumes many of the framework studied previously, and on the \textit{Sinkhorn divergence barycenters} introduced in \cite{Janati2020} (see also  \cite{Feydy2019}). 

\begin{dfn}[Doubly regularized Wasserstein barycenters]
Given a collection of weights $\blambda \in \Delta^n$, the \emph{doubly regularized  $W_p$-barycenter} of a family of measures $\bnu \in \cP_p(\R^d)^n$ with 
regularization parameters $\alpha,\beta \in \R_+$ is defined by 
\begin{equation}
\EBar_p^{\blambda}(\bnu \,|\, \alpha,\beta) := \argmin\limits_{\mu \in \cP_p(\R^d)} \sum_{i=1}^n \lambda_i \, \Eot_p(\nu_i,\mu \,|\, \alpha) + \beta\cdot \mathrm{KL}(\mu \,|\, \Leb)
\end{equation}
where $\Leb \in \cM_{\mathrm{loc}}(\R^d)_+$ denotes the Lebesgue measure.
\end{dfn}

\begin{dfn}[Sinkhorn divergence barycenter]
Given a collection of weights $\blambda \in \Delta^n$, the $p$-\emph{Sinkhorn divergence barycenter} of a collection of measures $\bnu \in \cP_p(\R^d)^n$ with regularization parameter $\alpha \in \R_+$ is defined by 
\begin{equation}
\SdBar_p^{\blambda}(\bnu \,|\, \alpha) := \argmin\limits_{\mu \in \cP_p(\R^d)} \sum_{i=1}^n \lambda_i \, \Eot_{p}(\nu_i,\mu \,|\, \alpha) - \frac{1}{2} \Eot_{p}(\mu,\mu \,|\, \alpha).
\end{equation}
\end{dfn}

In what follows, we prove that all the barycenter-type Wasserstein aggregators introduced above have the same breakdown point measured in the $W_p$-metric as the usual arithmetic mean. 

\begin{thm}[Breakdown point of double regularizd Wasserstein barycenters]
\label{thm:BreakdownWass1}
For each couple $(\blambda,\bnu) \in \Delta^n \times \cP_p(\R^d)^n$ and every $\alpha,\beta \in \R_+$, one has that 
\begin{equation*}
\bbeta^n_{W_p} \Big(\EBar_p^{\blambda}(\bnu \,|\, \alpha,\beta)  \Big) = \min \Big\{\lambda_i \,|\, i \in \llbracket 1,n \rrbracket, \; \lambda_i > 0\Big\}, 
\end{equation*}
and $\bbeta^{\infty}_{W_p}(\EBar_p(\bnu \,|\, \alpha,\beta)) = 0$. The same conclusions hold for $\SdBar_p^{\blambda}(\bnu \,|\,\alpha)$. 
\end{thm}

The proof of the above theorem is based on an intermediate first-order optimality condition for Wasserstein barycenters, involving the so-called duality mapping 
\begin{equation*}
j_p(z) := \, \|z\|^{p-2} z
\end{equation*}
defined for $p \in (1,+\infty)$. In the particular case in which $p=2$, we recover a known characterization of the space average  
\begin{equation*}
\E[\mu] := \int_{\R^d} x \, \mathrm{d} \mu(x)    
\end{equation*}
of the above aggregators, which was already found e.g. in \cite[Remark 2.6]{Carlier2021}.  

\begin{prop}[A first-order optimality condition]
\label{prop:CN1Bar}
For each $\bar{\mu}_{\alpha,\beta} \in \EBar_p^{\blambda}(\bnu \,|\, \alpha,\beta)$ and every choice of optimal plans $\bar{\gamma}_i \in \Gamma_o^{\alpha}(\bar{\mu}_{\alpha,\beta},\nu_i)$ with $i \in \llbracket 1,n \rrbracket$, one has that
\begin{equation}
\label{eq:BarycenterOptimality}
\sum_{i=1}^n \lambda_i \int_{\R^{2d}} j_p(x-y) \, \mathrm{d}\bar{\gamma}_i(x,y) = 0.
\end{equation}
When $p=2$, this entails in particular  
\begin{equation}
\label{eq:SpaceAverage}
\E [ \bar{\mu}_{\alpha,\beta} ] = \sum_{i=1}^n \lambda_i \E[\nu_i].
\end{equation}
Both the above formulas hold for $\bar{\mu}_{\alpha} \in \SdBar_p^{\blambda}(\bnu \, |\,\alpha)$ instead. 
\end{prop}

\begin{proof}
We will derive \eqref{eq:BarycenterOptimality} by comparing the costs between any $\bar{\mu}_{\alpha,\beta} \in \EBar_p^{\blambda}(\bnu \,|\, \alpha,\beta)$ and the translated competitor $\bar{\mu}^h_{\alpha,\beta} := (\Id+h)_{\sharp} \bar{\mu}_{\alpha,\beta}$ defined for some $h \in \R^d$.  If $\beta =0$, then the outer regularization term disappears from the cost, so suppose for now that $\beta > 0$. Hence necessarily $\bar{\mu}_{\alpha,\beta},\bar{\mu}_{\alpha,\beta}^h \ll \Leb$, so that by translation invariance of the Lebesgue measure, it holds that
\begin{equation}
\label{eq:ProofManipulability1}
\begin{aligned}
\mathrm{KL}(\bar{\mu}^h_{\alpha,\beta} \,|\, \Leb) & = \mathrm{KL}\Big((\Id+h)_{\sharp} \bar{\mu}_{\alpha,\beta} \,|\,(\Id+h)_{\sharp} \Leb \Big) \\
& = \int_{\R^{d}} \log \bigg( \frac{\mathrm{d}(\Id+h)_{\sharp} \bar{\mu}_{\alpha,\beta}}{\mathrm{d}(\Id+h)_{\sharp} \Leb} \bigg) \, \mathrm{d} (\Id+h)_{\sharp} \bar{\mu}_{\alpha,\beta} \\
& = \int_{\R^{d}} \log \bigg( \frac{\mathrm{d}\bar{\mu}_{\alpha,\beta}}{\mathrm{d} \Leb} \circ (\Id-h)\bigg) \, \mathrm{d} (\Id+h)_{\sharp} \bar{\mu}_{\alpha,\beta} \\
& = \mathrm{KL}(\bar{\mu}_{\alpha,\beta} \,|\, \Leb).
\end{aligned}
\end{equation}
Now, pick for each $i \in \llbracket 1,n \rrbracket$ some optimal plan $\bar{\gamma}_i \in \Gamma_o(\bar{\mu}_{\alpha,\beta},\nu_i)$ for the entropic transport problem, and let $\bar{\gamma}_i^h := (\pi^1+h,\pi^2)_{\sharp} \bar{\gamma}_i \in \Gamma(\bar{\mu}^h_{\alpha,\beta},\nu_i)$. Upon noticing that $\bar{\mu}_{\alpha,\beta}^h \otimes \nu_i = (\pi^1+h ,\pi^2)_{\sharp} (\bar{\mu}_{\alpha,\beta} \otimes \nu_i)$, it follows from similar computations as above that
\begin{equation}
\label{eq:ProofManipulability2}
\begin{aligned}
\mathrm{KL}(\bar{\gamma}_i^h \,|\, \mu_{\alpha,\beta}^h \otimes \nu_i) & = \int_{\R^{2d}} \log \bigg( \frac{\mathrm{d} \bar{\gamma}_i^h}{\mathrm{d} (\bar{\mu}_i^h \otimes \nu_i)} \bigg) \, \mathrm{d} \bar{\gamma}_i^h \\
& = \int_{\R^{2d}} \log \bigg( \frac{\mathrm{d} \bar{\gamma}_i}{\mathrm{d} (\bar{\mu}_i \otimes \nu_i)} \circ (\pi^1-h,\pi^2) \bigg) \, \mathrm{d} \bar{\gamma}_i^h \\
& = \mathrm{KL}(\bar{\gamma}_i \,|\, \mu_{\alpha,\beta} \otimes \nu_i). 
\end{aligned}
\end{equation}
In addition, it follows from standard superdifferentiability estimates for powers of the Euclidean norm (see e.g. \cite[Lemma 10.2.1]{AGS}) that 
\begin{equation}
\label{eq:ProofManipulability3}
\begin{aligned}
\int_{\R^{2d}} \|x-y\|^p \, \mathrm{d}( \bar{\gamma}_i^h - \bar{\gamma}_i)(x,y) & = \int_{\R^{2d}} \Big( \|x+h-y\|^p - \|x-y\|^p \Big)\, \mathrm{d} \bar{\gamma}_i)(x,y) \\
& \leq p \int_{\R^{2d}} \langle j_p(x-y) , h \rangle \, \mathrm{d} \bar{\gamma}_i(x,y) + \omega_p(h)
\end{aligned}
\end{equation}
with $\omega_p(h) = o(\|h\|)$ as $\|h\| \to 0^+$. Thence, upon combining \eqref{eq:ProofManipulability1}, \eqref{eq:ProofManipulability2} and \eqref{eq:ProofManipulability3} together with the optimality of $\bar{\mu}_{\alpha,\beta} \in \EBar_p^{\blambda}(\bnu \,|\, \alpha,\beta)$ and the arbitrariness of $h \in \R^d$, we get that 
\begin{equation}
\label{eq:FirstOrderOptProof}
\sum_{i=1}^n \lambda_i \int_{\R^{2d}} j_p(x-y) \, \mathrm{d} \bar{\gamma}_i(x,y) = 0
\end{equation}
for every choice of optimal plans $\bar{\gamma}_i \in \Gamma_o^{\alpha}(\bar{\mu}_{\alpha,\beta},\nu_i)$ with $i \in \llbracket 1,n \rrbracket$. In particular, since $j_2 = \Id$, we directly recover the mean value identity \eqref{eq:SpaceAverage} from the previous expression when $p=2$. Finally, the proof is identical in the case $\beta =0$, as we simply drop the outer entropic regularization. 

The case $\bar{\mu}_{\alpha} \in \SdBar_p^{\blambda}(\bnu \,|\,\alpha)$ is almost identical, with the exception of the entropic self-distance term. To handle this one, we simply notice that
\begin{equation*}
\Gamma(\bar{\mu}^h_{\alpha},\bar{\mu}^h_{\alpha}) = \Big\{(\pi^1+h,\pi^2+h)_{\sharp} \gamma \,\big|\, \gamma \in \Gamma(\bar{\mu}_{\alpha},\bar{\mu}_{\alpha}) \Big\} 
\end{equation*}
since $(\Id+h)$ is a bijection. Moreover, it can be shown through similar computations as \eqref{eq:ProofManipulability2} that 
\begin{equation*}
\mathrm{KL}(\bar{\gamma}_{\alpha} \,|\, \bar{\mu}_{\alpha} \otimes \bar{\mu}_{\alpha}) = \mathrm{KL} \Big( (\pi^1+h,\pi^2+h)_{\sharp} \bar{\gamma}_{\alpha} \,|\, \bar{\mu}^h_{\alpha} \otimes \bar{\mu}^h_{\alpha}\Big)
\end{equation*}
for every $\bar{\gamma} \in \Gamma(\bar{\mu}_{\alpha},\bar{\mu}_{\alpha})$. Thus, because $(\pi^1+h,\pi^2+h)$ is also an isometry, we have that  
\begin{equation*}
\int_{\R^{2d}} \|x-y\|^p \, \mathrm{d} \bar{\gamma}_{\alpha}(x,y) = \int_{\R^{2d}} \|x-y\|^p \, \mathrm{d} \Big( (\pi^1+h,\pi^2+h)_{\sharp}\bar{\gamma}_{\alpha} \Big)(x,y)    
\end{equation*}
which eventually yields
\begin{equation*}
\Eot_p(\bar{\mu}_{\alpha}^h,\bar{\mu}_{\alpha}^h \,|\, \alpha) = \Eot_p(\bar{\mu}_{\alpha},\bar{\mu}_{\alpha} \,|\, \alpha).
\end{equation*}
In conclusion, the self-entropic term is invariant under the pushforward by $(\Id+h)$, and from there we recover \eqref{eq:BarycenterOptimality} by repeating the arguments detailed above. 
\end{proof}

\begin{proof}[Proof of Theorem \ref{thm:BreakdownWass1}]
We begin by proving the statements for $\EBar_p^{\blambda}(\bnu\,|\,\alpha,\beta)$. Take some index $j \in \argmin\{\lambda_i \,|\, \lambda_i > 0\}$, pick any $\bsigma \in \cP_p(\R^d)^n$ so that $\sigma_i = \nu_i$ for $i \neq j$ and $\sigma_j = \delta_{y_j}$, and let $\tilde{\mu}_{\alpha,\beta}^j \in \EBar_p^{\blambda}(\bsigma \,|\,\alpha,\beta)$. Suppose now by contradiction that
\begin{equation*}
W_p(\bar{\mu}_{\alpha,\beta} ,\tilde{\mu}_{\alpha,\beta}^j) \leq M 
\end{equation*}
for some constant $M >0$ that is independent of $y_j \in \R^d$, which  by \eqref{eq:CrudeBound} in Lemma \ref{lem:PseudoTriangle} entails in particular that 
\begin{equation}
\label{eq:MomentBound}
\cM_p(\tilde{\mu}_{\alpha,\beta}^j) \leq C_M \qquad \text{and} \qquad \Eot_p(\bar{\mu}_{\alpha,\beta}^j,\tilde{\mu}_{\alpha,\beta}^j \,|\, \alpha) \leq C_M'
\end{equation}
for some other constants $C_M,C_M' > 0$ depending only on $M >0$ and $\bnu \in \cP_p(\R^d)^n$. From \eqref{eq:FirstOrderOptProof} combined with the (pseudo) triangle inequality satisfied by the entropic optimal transport problem stated in Lemma \ref{lem:PseudoTriangle}, this implies that
\begin{equation}
\label{eq:BlowupProof0}
\begin{aligned}
\bigg\| \, \lambda_j \cdot \int_{\R^{2d}} j_p(y-x) \, \mathrm{d} \tilde{\gamma}_j(x,y) \, \bigg\| & = \bigg\| \, \sum_{i \neq j} \lambda_i \cdot \int_{\R^{2d}} j_p(x-y) \, \mathrm{d} \tilde{\gamma}_i(x,y) \, \bigg\| \\
& \leq \sum_{i \neq j} \lambda_i \cdot \Eot_p(\tilde{\mu}_{\alpha,\beta}^j,\nu_i \,|\, \alpha)^{\nicefrac{1}{q}} \\
& \leq 2^{\nicefrac{(p-1)}{q}} \sum_{i \neq j} \lambda_i \cdot \Big( \Eot_p(\tilde{\mu}_{\alpha,\beta}^j,\bar{\mu}_{\alpha,\beta} \,|\, \alpha) + \Eot_p(\bar{\mu}_{\alpha,\beta},\nu_i \,|\, \alpha)\Big)^{\nicefrac{1}{q}} \\
& \leq C_1
\end{aligned}
\end{equation}
for every choice of optimal plans $\tilde{\gamma}_i \in \Gamma_o^{\alpha}(\tilde{\mu}_{\alpha,\beta}^j,\nu_i)$ with $i \in \llbracket 1,n \rrbracket$ and some constant $C_1 > 0$ depending on $M$, $p$ and $\bnu \in \cP_p(\R^d)^n$. On the other hand, upon noting that $\Gamma(\tilde{\mu}_{\alpha,\beta}^j,\delta_{y_j}) = \{\tilde{\mu}_{\alpha,\beta}^j \otimes \delta_{y_j}\}$, it holds that 
\begin{equation}
\label{eq:BlowUpProof1}
\begin{aligned}
\bigg\langle \frac{y_j}{\|y_j\|} \, , \, \int_{\R^{2d}} j_p(y-x) \, \mathrm{d} \tilde{\gamma}_j(x,y) \bigg\rangle & = \int_{\{\|x\| \leq \|y_j\|/2 \}} \|y_j-x\|^{p-2} \Big\langle \frac{y_j}{\|y_j\|} , y_j - x \Big\rangle \, \mathrm{d} \tilde{\mu}_{\alpha,\beta}^j(x) \\
& \hspace{0.45cm} + \int_{\{\|x\| > \|y_j\|/2 \}} \|y_j-x\|^{p-2} \Big\langle \frac{y_j}{\|y_j\|} , y_j - x \Big\rangle \, \mathrm{d} \tilde{\mu}_{\alpha,\beta}^j(x) \\
& \geq \int_{\{\|x\| \leq \|y_j\|/2 \}} \|y_j-x\|^{p-2} \Big( \|y_j\| - \|x\| \Big) \, \mathrm{d} \tilde{\mu}_{\alpha,\beta}^j(x) \\
& \hspace{0.45cm} - \int_{\{\|x\| > \|y_j\|/2 \}} \|y_j-x\|^{p-1} \, \mathrm{d} \tilde{\mu}_{\alpha,\beta}^j(x) 
\end{aligned}
\end{equation}
At this point, the second integral term in the right-hand side of the above expression may be further estimated from below using Chebyshev's and Hölder's inequalities, by noting that 
\begin{equation}
\label{eq:BlowUpProof2}
\begin{aligned}
\int_{\{\|x\| > \|y_j\|/2 \}} \|y_j-x\|^{p-1} \, \mathrm{d} \tilde{\mu}_{\alpha,\beta}^j(x) & \leq \tilde{\mu}_{\alpha,\beta}^j \Big( \Big\{ x \in \R^d \, \big| \, \|x\| > \tfrac{\|y_j\|}{2} \Big\}\Big)^{\nicefrac{1}{p}} \cdot  W_p(\tilde{\mu}_{\alpha,\beta}^j,\delta_{y_j})^{\nicefrac{p}{q}} \\
& \leq \frac{2 C_M^p}{\|y_j\|} \cdot \Big( \|y_j\| + C_M \Big)^{\nicefrac{p}{q}} \\
& \leq C_2 \cdot \|y_j\|^{p-2} \Big( 1 + \tfrac{C_M}{\|y_j\|} \Big)^{\nicefrac{p}{q}}
%
%
\end{aligned}
\end{equation}
for some constant $C_2 > 0$ depending again only on $M$, $p$ and $\bnu \in \cP_p(\R^d)^n$. Moreover, the first integral in the last estimate of \eqref{eq:BlowUpProof1} can also be lower-bounded by using Chebyshev's inequality, the estimates of \eqref{eq:MomentBound} as 
\begin{equation}
\label{eq:BlowUpProof3}
\begin{aligned}
\int_{\{\|x\| \leq \|y_j\|/2 \}} \|y_j-x\|^{p-2} \Big( \|y_j\| - \|x\| \Big) \, \mathrm{d} \tilde{\mu}_{\alpha,\beta}^j(x) & \geq C_p \|y_j\|^{p-1}  \tilde{\mu}_{\alpha,\beta}^j \Big( \Big\{ x \in \R^d \, \big| \, \|x\| \leq \tfrac{\|y_j\|}{2} \Big\}\Big) \\
& = C_p \|y_j\|^{p-1} \bigg( 1 -  \tilde{\mu}_{\alpha,\beta}^j \Big( \Big\{ x \in \R^d \, \big| \, \|x\| > \tfrac{\|y_j\|}{2} \Big\}\Big) \bigg) \\
& \geq C_p \|y_j\|^{p-1} \Big( 1 -  \big( \tfrac{2 C_M}{\|yj\|}  \big)^p \Big)
\end{aligned}
\end{equation}
with $C_p := \tfrac{1}{2} \min\{(1/2)^{p-2},(3/2)^{p-2}\}$, where we used the monotonicity of the map $t \mapsto t^{p-2}$ for $p \in (1,+\infty)$, and the fact that $\|y_j\|/2 \leq \|x-y_j\| \leq 3\|y_j\|/2$ on the above integration domain. Whence, merging \eqref{eq:BlowUpProof1} and \eqref{eq:BlowUpProof2} in \eqref{eq:BlowUpProof3} and letting $\|y_j\| \to +\infty$ eventually violates the uniform upper-bound derived in \eqref{eq:BlowupProof0}.
\end{proof}


Another way of aggregating measures proposed in \cite{Alvarez2016,Werenski2024} is to replace the Wasserstein distances by their linearizations around a given reference measure $\rho \in \cP_p^r(\R^d)$.

\begin{dfn}[Linearized Wasserstein barycenters]
Given a collection of weights $\blambda \in \Delta^n$, the $p$-\emph{linearized Wasserstein barycenter} around $\rho \in \cP_p^{\mathrm{ac}}(\R^d)$ of a family of measures $\bnu \in \cP_p(\R^d)^n$ is defined by 
\begin{equation}
\LWBar_p^{\blambda}(\bnu \,|\,\rho) = \argmin\limits_{\mu \in \cP_p(\R^d)} \sum_{i=1}^n \lambda_i \cdot W_{p,\rho}^p(\nu_i,\mu)
\end{equation}
\end{dfn}

We end this section by showing that linearized 2-Wasserstein barycenters also have the same breakdown point. While we conjecture that this statement remains true for $p \in (1,+\infty)$, we were not able to establish a clean first-order optimality condition akin to Proposition \ref{prop:CN1Bar} supporting it. 

\begin{thm}[Breakdown point of linearized 2-Wasserstein barycenters]
\label{thm:BreakdownWass2}
For each couple $(\blambda,\bnu) \in \Delta^n \times \cP_2(\R^d)^n$ and every $\rho \in \cP_2^{\mathrm{ac}}(\R^d)$, one has that 
\begin{equation*}
\bbeta^n_{W_{2,\rho}} \Big(\LWBar_2^{\blambda}(\bnu \,|\, \rho)  \Big) = \bbeta^n_{W_2} \Big(\LWBar_2^{\blambda}(\bnu \,|\, \rho)  \Big) = \min\Big\{\lambda_i \,|\, i \in \llbracket 1,n \rrbracket, \; \lambda_i > 0\Big\}, 
\end{equation*}
and likewise $\bbeta^{\infty}_{W_{2,\rho}}(\LWBar_2( \bnu \,|\, \rho)) = \bbeta^{\infty}_{W_2}(\LWBar_2( \bnu \,|\, \rho))  = 0$. 
\end{thm}

\begin{proof}
The demonstration of this result hinges upon the fact that every $\bar{\mu}_{\rho} \in \LWBar_2^{\blambda}(\bnu \,|\,\rho)$ is of the form $\bar{\mu}_{\rho} = \bar{T}_{\rho}(\bnu)_{\sharp} \rho$, with  
\begin{equation}
\label{eq:LinBarIdentity}
\bar{T}_{\rho}(\bnu) = \argmin \bigg\{ \sum_{i=1}^n \lambda_i \cdot \|T - T_{\rho}^{\nu_i}\|_{L^2(\rho)}^2 \, \Big|\, T \in L^2(\R^d,\R^d;\rho) \cap \nabla \cC(\R^d) \bigg\}, 
\end{equation}
Therein $\nabla \cC(\R^d)$ stands for the (convex) set of gradients of convex functions, and as it follows e.g. from \cite[Corollary 5.23]{Villani2} that Brenier maps are closed under strong $L^2$-convergence, the above minimization problem is well-posed and admits a unique solution of the form 
\begin{equation*}
\bar{T}_{\rho}(\bnu) = \sum_{i=1}^n \lambda_i \cdot T_{\rho}^{\nu_i}.
\end{equation*}
This identity, which was already derived in \cite{Werenski2024}, stems from the observation that the map $\bar{T}_{\rho}(\bnu)$ defined above solves the relaxed problem 
\begin{equation*}
\min \bigg\{ \sum_{i=1}^n \lambda_i \cdot \|T - T_{\rho}^{\nu_i}\|_{L^2(\rho)}^2 \, \Big|\, T \in L^2(\R^d,\R^d;\rho) \bigg\}
\end{equation*}
while remaining admissible for \eqref{eq:LinBarIdentity} by convexity of the cone $\nabla \cC(\R^d)$. With this characterization at hand, fix any $j \in \argmin_{i\in\llbracket 1,n \rrbracket}\{\lambda_i \,|\, \lambda_i > 0\}$, take $\bsigma \in \cP_2(\R^d)^n$ be such that 
$\sigma_i = \nu_i$ for $i \neq j$ and $\sigma_j = \delta_{y_j}$, and let $\tilde{\mu}_{\rho} = \LWBar_2^{\blambda}(\bsigma \,|\, \rho)$. Then, following Theorem \ref{thm:Brenier}, it holds that 
\begin{equation}
\label{eq:LinearizedBlowup}
\begin{aligned}
W_{2,\rho}(\bar{\mu}_{\rho},\tilde{\mu}_{\rho}^j) & = \bigg\| \sum_{i=1}^n \lambda_i \cdot (T_{\rho}^{\nu_i} - T_{\rho}^{\sigma_i}) \bigg\|_{L^2(\rho)} \\
& = \lambda_j \cdot \| T_{\rho}^{\nu_j} - y_j \|_{L^2(\rho)} ~\underset{\|y_j\| \to +\infty}{\longrightarrow}~ +\infty. 
\end{aligned}
\end{equation}
This shows that $\bbeta^n_{W_{2,\rho}}(\LWBar_2^{\blambda}(\bnu)) = \min\{ \lambda_i \,|\, \lambda_i > 0\}$, which by definition of the asymptotic breakdown point directly implies that $\bbeta^{\infty}_{W_{2,\rho}}(\LWBar_2^{\blambda}(\bnu)) = 0$. Concerning the breakdown point measured in the standard $W_2$-metric, it suffices to observe that 
\begin{equation*}
\begin{aligned}
\E[\bar{\mu}_{\rho}] & = \int_{\R^d} x \, \mathrm{d}\bigg( \mathsmaller{\sum}\limits_{i=1}^n \lambda_i \cdot T_{\rho}^{\nu_i} \bigg)_{\raisebox{4pt}{{\scriptsize $\sharp$}}} \rho(x) \\
& = \sum_{i=1}^n \lambda_i \cdot \int_{\R^d} T_{\rho}^{\nu_i}(x) \, \mathrm{d} \rho(x) =  \sum_{i=1}^n \lambda_i \cdot \E[\nu_i]
\end{aligned}
\end{equation*}
which combined with the basic estimate 
\begin{equation*}
W_2(\bar{\mu}_{\rho},\tilde{\mu}_{\rho}^j) \geq \lambda_j \cdot \| \E[\bar{\mu}_{\rho}] - \E[\tilde{\mu}_{\rho}^j] \| = \lambda_j \cdot \| \E[\nu_j] - y_j \| 
\end{equation*}
provides the result by the same contradiction argument as above. 
\end{proof}

\begin{rem}[Absence of exact manipulability for Wasserstein barycenters]
The breakdown point is inherently a qualitative measure of the robustness of an aggregator, and as such does not say anything about the collection of states which are effectively reachable by corrupting a sufficient fraction of input data. To illustrate this fact, it is well-known that the arithmetic mean in Euclidean spaces is exactly manipulable by one adversary with full-information, in the sense that given any faithful input $\bx \in \cH^n$ and target $\bar{x}' \in \cH$, picking an index  $j \in \llbracket 1,n \rrbracket$, the adverse input $\by \in \cH^n$ defined by $y_i = x_i$ for every $i \neq j$ and 
\begin{equation*}
y_j := x_j + \frac{1}{\lambda_j} \Big( \bar{x}' - \Mean^{\blambda}(\bx) \Big)
\end{equation*}
is such that $\Mean^{\blambda}(\by) = \bar{x}'$.  

Such an arbitrary manipulability fails in general for Wasserstein barycenters, both classical and linearized. Indeed, it follows e.g. from \cite[Theorem 5.1]{Agueh2011} that if at least one of the elements in $\bnu \in \cP_2(\R^d)^n$ is absolutely continuous with a bounded density, so is $\WBar_2^{\lambda}(\bnu) \in \cP_2(\R^d)$. Therefore, it is not possible to reach an arbitrary target measure $\bar{\mu}' \in \cP_2(\R^d)$ by modifying just one or a few inputs. Such an exact manipulability seems however to partially restricted to some particular classes of measures. For instance, if all the inputs are Dirac masses, then the Wasserstein barycenter reduces to a Euclidean one and every target Dirac measure is reachable. Similarly, if all the inputs $\bnu \in \cP_2(\R^d)^n$ are isotropic Gaussians, namely $\nu_i = \cN(m_i,\sigma_i \cdot \Id)$ for each $i \in \llbracket 1,n \rrbracket$, then it follows from characterization provided in \cite[Theorem 6.1]{Agueh2011} of their barycenters that $\WBar_2^{\blambda}(\bnu) = \cN(\bar{m},\bar{\sigma} \cdot \Id)$, with 
\begin{equation*}
\bar{m} = \sum_{i=1}^n \lambda_i \cdot m_i \qquad \text{and} \qquad \bar{\sigma} = \bigg( \sum_{i=1}^n \lambda_i \cdot \sigma_i^{\nicefrac{1}{2}} \bigg)^2
\end{equation*}
Thus, upon perturbing one input, it is possible to reach any other isotropic Gaussian $\cN(\bar{m}',\bar{\sigma}' \cdot \Id)$ with arbitrary mean $\bar{m}' \in \R^d$, but whose covariance must satisfy $\bar{\sigma}' > (\sum_{i\neq j} \lambda_i \cdot \sigma_i^{\nicefrac{1}{2}})^2$.
\end{rem}


\subsection{Hölder and Lipschitz resilient aggregators in metric spaces}
\label{section:Resilient}

In this section, we propose a general notion of quantitative resilience against adversaries for weighted aggregators defined over metric spaces, inspired by the notion of Lipschitz resilience recently put forth in \cite{Allouah2024}. 

\begin{dfn}[Hölder and Lipschitz resilience]
\label{def:Resilience}
Given $\blambda \in \Delta^n$, we say that single-valued weighted aggregator $\scF^{\blambda} : \cX^n \to \cX$ is \emph{Hölder resilient} with constant $L > 0$ and exponent $\alpha \in (0,1]$ provided 
\begin{equation*}
\dist \Big(\scF^{\blambda}(\bx),\scF^{\blambda}(\by) \Big) \leq L \cdot \CmD_{\blambda}(\bx,\by)^\alpha =  L \cdot \bigg(\sum_{i \in \msupp{\bx}{\by}} \lambda_i\bigg)^\alpha
\end{equation*}
for all $\bx,\by \in \cX^n$. In words, the maximum influence of a subset $\msupp{\bx}{\by}$ of mismatches is proportional to a power of their cumulated weights. In particular, we say that $\scF^{\blambda}$ is \textit{Lipschitz resilient} if it is Hölder resilient with exponent $\alpha = 1$.
\end{dfn}

\begin{rem}[On the breakdown point of H\"older resilient aggregators]
It is straightforward to see that in a metric space with infinite diameter, the breakdown point of a Hölder resilient aggregator is always $+\infty$, since $\dist_H(\scF^{\blambda}(\bx), \scF^{\blambda}(\by))$ is controlled by a quantity that does not depend on the inputs and therefore cannot go to infinity. 
The price to pay for this very strong robustness is that Hölder resilient aggregator cannot be exact, in the sense that one typically expects that $\scF^{\blambda}(\bx) \neq x$ for $\bx := (x,x,\ldots,x) \in \cX^n$, as the aggregator maps the whole of $\cX^n$ into a set of fixed diameter. Such an aggregator similarly cannot be equivariant under the action of an unbounded group of isometries (see Proposition~\ref{prop:breakdown} above). While this may seem restrictive and rather conservative, Lipschitz resilient aggregators are are currently being used for robust content recommendations on the platform \href{https://tournesol.app}{Tournesol.app} through a Bayesian generalization of Bradley-Terry models, see~\cite{Fageot2024} for details.
\end{rem}

In the following theorem, we exhibit a fairly large class of H\"older resilient aggregators, which we shall particularize to Wasserstein spaces in Section \ref{section:QrWm} below. 

\begin{thm}[A general family of H\"older resilient aggregators]
\label{thm:Resilience}
Let $\blambda \in \Delta^n$ and $\scF^{\blambda} : \cX^n \to \cX$ be a weighted aggregator defined through the minimization problem 
\begin{equation*}
\scF^{\blambda}(\bx) := \argmin_{y \in \cX} \cJ_L^{\blambda}(\bx,y)
\end{equation*}
where the cost function $\cJ_L^{\blambda} : \cX^n \times \cX \to \R_+$ is of the form
\begin{equation*}
\cJ_L^{\blambda}(\bx,y) := \psi_L(y) + \sum_{i=1}^n \lambda_i \cdot \dsf_i(x_i,y)   
\end{equation*}
and satisfies the following assumptions. 

\begin{enumerate}
\item[$(i)$] For each $i \in \llbracket 1,n \rrbracket$, the map $\dsf_i : \cX^2 \longrightarrow \R_+$, is a pseudo-metric satisfying 
\begin{equation*}
\dsf_i(x,y) \leq \dist(x,y)
\end{equation*}  
for all $x,y \in \cX$.
\item[$(ii)$] The map $\psi_L : \cX \to \R_+$ is such that $\cJ_L^{\blambda}$ satisfies the generalized uniform convexity estimate
\begin{equation*}
\cJ_L^{\blambda}(\bx,y) \geq \cJ_L^{\blambda}( \bx,\scF^{\blambda}(\bx)) + \frac{1}{L} \cdot \dist\Big( \scF^{\blambda}(\bx) , y \Big)^{1+\frac{1}{\alpha}}. 
\end{equation*}
for some $\alpha \in (0,1]$. 
\end{enumerate}
Then $\scF^{\blambda} : \cX^n \to \cX$ is $\alpha$-H\"older resilient with 
\begin{equation*}
\dist\left(\vphantom{\sum}\scF^{\blambda}(\bx), \scF^{\blambda}(\by) \right) \leq \Big( L \cdot \CmD_{\blambda}(\bx,\by) \Big)^{\alpha}.
\end{equation*}
Moreover, denoting by $\bar{x}_L := \argmin \cJ_L^{\blambda}(\bx,y)$, it holds that 
\begin{equation}
\label{eq:DataFidelty}
\dist(\bar{x}_L,x')^{1+\frac{1}{\alpha}} \leq L \cdot \Big( \psi_L(x') - \psi_L(\bar{x}_L) \Big)
\end{equation}
for each $x' \in \argmin\limits_{y \in \cX} \sum\limits_{i=1}^n \lambda_i \cdot \dsf_i(x_i,y)$.
\end{thm}

\begin{proof}
We begin by observing that the statement is obviously true if $\dist(\scF^{\blambda}(\bx), \scF^{\blambda}(\by)) = 0$, so in the sequel we posit that $\dist(\scF^{\blambda}(\bx), \scF^{\blambda}(\by)) > 0$. By the generalized uniform convexity of item $(ii)$, it holds for any $(\bz,y) \in \cX^n \times \cX$ that 
\begin{align*}
\dist \Big( \scF^{\blambda}(\bz),y \Big)^{1+\frac{1}{\alpha}} & \leq L \cdot \Big( \cJ_L^{\blambda}(\bz,y) - \cJ_L^{\blambda}(\bz,\scF^{\blambda}(\bz)) \Big) 
\end{align*}
so that, alternately picking $(\bz,y) := (\bx,\scF^{\blambda}(\by))$ and $(\bz,y) := (\by,\scF^{\blambda}(\bx))$ and summing the resulting two inequalities, we obtain
\begin{subequations}
\begin{align}
2 \cdot \dist\left(\vphantom{\sum}\scF^{\blambda}(\bx), \scF^{\blambda}(\by) \right)^{1+\frac{1}{\alpha}} & \leq L \cdot \Big( \cJ_L^{\blambda}(\by,\scF^{\blambda}(\bx)) - \cJ_L^{\blambda}(\by,\scF^{\blambda}(\by)) \Big) \nonumber \\
& \hspace{1.15cm} + L \cdot \Big( \cJ_L^{\blambda}(\bx,\scF^{\blambda}(\by)) - \cJ_L^{\blambda}(\bx,\scF^{\blambda}(\bx)) \Big) \nonumber \\
& = L \cdot \Big( \cJ_L^{\blambda}(\by,\scF^{\blambda}(\bx)) - \cJ_L^{\blambda}(\bx,\scF^{\blambda}(\bx)) \Big) \nonumber \\
& \hspace{1.15cm} + L \cdot \Big( \cJ_L^{\blambda}(\bx,\scF^{\blambda}(\by))- \cJ_L^{\blambda}(\by,\scF^{\blambda}(\by)) \Big) \nonumber \\
& = L \cdot  \sum_{i \in \msupp{\bx}{\by}} \lambda_i \cdot \Big( \dsf_i(\scF^{\blambda}(\bx), y_i) - \dsf_i(\scF^{\blambda}(\bx), x_i) \Big) \nonumber \\
& \hspace{1.15cm} + L \cdot  \sum_{i\in \msupp{\bx}{\by}} \lambda_i \cdot \Big( \dsf_i(\scF^{\blambda}(\by), x_i) - \dsf_i(\scF^{\blambda}(\by), y_i) \Big) \nonumber \\
& = L \cdot  \sum_{i\in\msupp{\bx}{\by}} \lambda_i \cdot \Big( \dsf_i(\scF^{\blambda}(\bx),y_i) - \dsf_i(\scF^{\blambda}(\by),y_i) \Big) \nonumber \\
& \hspace{1.15cm} + L \cdot \sum_{i\in\msupp{\bx}{\by}} \lambda_i \cdot \Big(  \dsf_i(\scF^{\blambda}(\by),x_i) - \dsf_i(\scF^{\blambda}(\bx),x_i) \Big) \nonumber \\
& \leq L \cdot \sum_{i\in\msupp{\bx}{\by}} \lambda_i \cdot \Big( \dsf_i(\scF^{\blambda}(\bx), \scF^{\blambda}(\by)) + \dsf_i(\scF^{\blambda}(\bx), \scF^{\blambda}(\by)) \Big) \label{step:triangle} \\
& \leq 2L \cdot \sum_{i\in\msupp{\bx}{\by}} \lambda_i \cdot \dist\left(\vphantom{\sum}\scF^{\blambda}(\bx), \scF^{\blambda}(\by) \right) \label{step:domin}
\end{align}
\end{subequations}
where~\eqref{step:triangle} stems from the triangle inequalities
\begin{equation*}
\left\{ 
\begin{aligned}
\dsf_i(\scF^{\blambda}(\bx),y_i) & \leq \dsf_i(\scF^{\blambda}(\bx),\scF^{\blambda}(\by)) + \dsf_i(\scF^{\blambda}(\by),y_i), \\
\dsf_i(\scF^{\blambda}(\by),x_i) & \leq \dsf_i(\scF^{\blambda}(\bx),\scF^{\blambda}(\by)) + \dsf_i(\scF^{\blambda}(\bx),x_i),
\end{aligned}
\right.
\end{equation*}
and~\eqref{step:domin} from the domination assumption made in $(i)$ on each pseudo-metric $\dsf_i$. As we assumed that $\dist(\scF^{\blambda}(\bx), \scF^{\blambda}(\by)) > 0$, we may divide both sides by the latter quantity and take the power $\alpha \in (0,1]$ to conclude. The estimate \eqref{eq:DataFidelty} may then be recovered through similar computations. 
\end{proof}

\begin{rem}[On the assumptions of Theorem \ref{thm:Resilience}]
Condition $(ii)$ in Proposition~\ref{thm:Resilience} is called a generalized uniform convexity property because, when $\alpha = 1$, it can be obtained e.g. as a direct consequence of the metric $\nicefrac{2}{L}$-strong convexity of $\psi_L : \cX \to \R_+$ along curves for which all the $\dsf_i$ are also convex.
\end{rem}

A direct consequence of Theorem~\ref{thm:Resilience} is the Lipschitz resilience of the quadratically regularized medians we introduced in the Introduction, and whose definition is recalled here.

\begin{dfn}[Quadratically regularized geometric median] 
\label{dfn:qrgm}
Given a regularization parameter $L>0$ and a collection of weights $\blambda \in (\R_+^{\star})^n$, the \emph{quadratically regularized geometric median} with anchor $x_0 \in \cH$ of a family of points $\bx = (x_1,\ldots,x_n) \in \cH^n$ is defined by 
\begin{equation} \label{eq:qrgm}
\QrGm^{\blambda}_L(\bx \,|\,x_0) := \argmin_{y \in \cH} \frac{1}{L} \cdot \|y - x_0\|^2_\cH + \sum_{i=1}^n \lambda_i \cdot \|y - x_i\|_{\cX}.
\end{equation}
\end{dfn}

\begin{rem}[A Bayesian perspective on quadratic regularization] \label{rem:bayes}
Although quadratic regularization may seem sort of arbitrary, it finds an elegant interpretation in terms of Bayesian statistics. Indeed, suppose that $(\cX,\dist) = (\R^d,\|\cdot\|)$ and the inputs $\bx =(x_1,\ldots,x_n) \subset (\R^d)^n$ are independently sampled from some Laplace distribution with parameter $\theta \in \R^d$, whose law is 
\begin{equation*}
\nu_\theta = \frac{1}{Z}\exp\big(\hspace{-0.05cm} -\|\cdot-\,\theta\|\, \big)  \cdot \Leb.
\end{equation*}
Then the Maximum Likelihood Estimator for the parameter $\theta \in \R^d$ maximizes the likelihood 
\begin{equation*}
\cL(\bx \,|\, \theta) :=  \frac{1}{Z^n} \prod_{i=1}^n \exp \big( \hspace{-0.05cm} -\|x_i - \theta\| \,\big) = \frac{1}{Z^n} \exp\bigg( \hspace{-0.075cm} -\sum_{i=1}^n \|x_i-\theta\| \, \bigg). 
\end{equation*}
Expressed in terms of the log-likelihood, this equivalently means that $\hat{\theta} = \GMed^{\bu_n}(\bx)$. Now, we posit the knowledge of a Bayesian prior over parameters, whose distribution is an isotropic Gaussian 
\begin{equation*}
\xi_0 := (\pi L)^{-\nicefrac{d}{2}} \exp\left(-\frac{\|\theta-x_0\|^2}{L} \right) \cdot \Leb.    
\end{equation*}
with mean $x_0 \in \R^d$ and variance $\nicefrac{L}{2} \in \R_+^{\star}$. Then, by virtue of Bayes law, we can derive the posterior distribution over $\theta \in \R^d$ as
\begin{equation*}
\xi(\theta \,|\, \bx) = \frac{\cL(\bx \,|\, \theta) \; \xi_0(\theta)}{\Pb(\bx)}
\end{equation*}
where the probability $\Pb(\bx)$ of observing a given input sample is independent from the parameter distribution. Hence, taking the log-posterior distribution yields
\begin{equation*}
\begin{aligned}
\log\xi(\theta \,|\, \bx) & = \log\cL(\bx \,|\,\theta) + \log \xi(\theta) - \log(\Pb(\bx)) \\
& = C - \frac{1}{L} \|\theta - x_0\|^2 - \sum_{i=1}^n \|\theta-x_i\|   
\end{aligned}
\end{equation*}
where $C = -\frac{d}{2}\log(\pi L) - n \log Z - \log(\Pb(\bx))$ is independent of $\theta \in \R^d$, so that the Maximum A Posteriori estimator is exactly $\Tilde{\theta} = \QrGm^{\one_n}_L(\bx \,|\,x_0)$. This gives a natural interpretation for the regularization parameter $L > 0$ as the confidence placed on the prior.
\end{rem}

Lipschitz resilient aggregators are hence closely related to bayesian statistical estimation. However, the above example is restricted to parametric estimation, where strong assumptions are made on the data distribution (here identically distributed according to a Laplace law for the data, with a Gaussian prior on its mean). Non-parametric estimation, e.g. with kernels or histograms, extends this to i.i.d. data with unknown common distribution. An interesting question relates to the extension of this framework to estimating agreement between non identically distributed variables.


\begin{ex}[$L^2$-regularized medians in $L^p$-spaces]
Take some $\rho \in \cP(\R^d)$ and fix an element $p \in [1,2]$. Then, the spaces defined by 
\begin{equation*}
\cH := L^2(\R^d,\R^d;\rho) \subset L^p(\R^d,\R^d;\rho) := \cX
\end{equation*}
are such that $\|\cdot\|_{\cX} \leq \|\cdot\|_{\cH}$, and one may thus define the quadratically regularized median 
\begin{equation*}
\QrGm^{\blambda}_{p,L}(\bof \,|\,f_0) := \argmin_{g \in L^2(\nu)} \frac{1}{L} \cdot \|g - f_0\|_{L^2(\rho)}^2 + \sum_{i=1}^n \lambda_i \cdot \|g - f_i\|_{L^p(\rho)}. 
\end{equation*}
\end{ex}

\begin{rem}[On the lack of monotony of $\QrGm^{\blambda}_L(\bx \,|\,x_0)$ restricted to the simplex]
\label{rem:NoMonotony}
As alluded to in Section \ref{section:Preli}, although quadratically regularized medians are aggregators in the sense of Definition \ref{def:Aggregators}, their restriction to the simplex $\Delta^n$ may fail to be monotone. Indeed, if we set $\blambda'_t := (\blambda+t e_j)/(1+t)$ for some $j \in \llbracket 1,n \rrbracket$ with $t > 0$ -- which amounts to increasing the weight $\lambda_j$ relative to the others and projecting back on the simplex --, then the aggregation problem can be rewritten as 
\begin{equation*}
\QrGm^{\blambda'_t}_L(\bx \,|\,x_0) = \argmin_{y \in \cH} \frac{1+t}{L} \cdot \|y - x_0\|^2_\cH + \sum_{i=1}^n \lambda_i \cdot \|y - x_i\|_{\cX} + t \cdot \|y - x_j\|_{\cX}.
\end{equation*}
It is then easy to see that the weight of the prior $x_0$ has also been increased, and since the latter is penalized quadratically, it may very well pull the aggregator away from the data point $x_j$. 
\end{rem}

In the next theorem, we show that quadratically regularized medians are Lipschitz resilient. 

\begin{thm}[Lipschitz resilience of quadratically regularized medians]
\label{thm:LipResGrGm}
For every $\blambda \in \Delta^n$, the aggregator $\QrGm^{\blambda}_L : \cH^n \to \cH$ is $L$-Lipschitz resilient.
\end{thm}

\begin{proof}
This is a direct consequence of Proposition~\ref{thm:Resilience} applied to the case $\dsf_i = \|\cdot\|_{\cX}$ and $\psi_L := \|\cdot - x_0\|_{\cH}^2/L$. Then, item $(i)$ of the proposition is clearly satisfied, whereas the cost function 
\begin{equation*}
\cJ_L^{\blambda}(\bx,y) = \frac{1}{L} \cdot \|y - x_0\|^2_\cH + \sum_{i=1}^n \lambda_i \, \|y - x_i\|_{\cX} 
\end{equation*}
is $\nicefrac{2}{L}$-strongly convex in $\|\cdot\|_{\cH}$ under our assumptions, which classically implies that the uniform convexity estimate of item $(ii)$ in the proposition holds with $\alpha =1$. 
\end{proof}

Another desirable feature of so above family of resilient aggregators is that, although biased towards the prior, they always remain closer to the unregularized one.

\begin{prop}[Fidelity to data]
It holds that 
\begin{equation*}
\big\| \QrGm_L^{\blambda}(\bx \,|\,x_0) - \GMed(\bx) \big\|_{\cH} \leq\, \big\| x_0 - \GMed(\bx) \big\|_{\cH}
\end{equation*}
for every $\bx \in \cH^n$. 
\end{prop}

\begin{proof}
This is an immediate consequence of \eqref{eq:DataFidelty} in Theorem \ref{thm:Resilience}. 
\end{proof}

We end this section on resilient aggregators by showing that, in the particular case in which $\cX = \cH$ is a Hilbert space, then both the geometric median and its quadratic regularization lie in the convex hull of the input data. 

\begin{lem}[Representer theorem for general medians in Hilbert spaces]
\label{lem:Hull}
\begin{subequations}
For every collection of weights $\blambda \in\Delta^n$, any multi-index $\ba \in [1,+\infty)^n$ and each family of points $\bx \in \cH^n$, it holds that
\begin{equation} \label{eq:convHull}
\argmin_{y\in\cH} \sum_{i=1}^n \lambda_i \cdot \|y - x_i\|_\cH^{a_i} \subset \conv(\bx).
\end{equation}
In particular, we have that
\begin{equation} \label{eq:convGMed}
\GMed^{\blambda}(\bx) \subset  \conv(\bx) \qquad \text{and} \qquad \QrGm_L^{\blambda}(\bx\,|\,x_0) \in \conv(\bx,x_0). 
\end{equation}
\end{subequations}
\end{lem}

\begin{proof}
Denote by 
\begin{equation*}
\cJ_L^{\blambda}(\bx,y \,|\,\ba) := \sum_{i=1}^n \lambda_i \cdot \|y - x_i\|_\cH^{a_i}
\end{equation*}
and let $\bar{x} \in \argmin\{\cJ_L^{\blambda}(\bx,y \,|\,\bp) \,|\, y \in \cH\}$. Observe first that if $\bar{x} = x_i$ for some $i \in \llbracket 1,n \rrbracket$, then the claim is automatically true. Suppose thus that $\bar{x} \notin \bx$, and assume by contradiction that $\bar{x} \notin \conv(\bx)$. Then, by the separating hyperplane theorem (see e.g. \cite[Theorem 1.7]{Brezis}), there exists a $p \in \cH$ such that $\langle \bar{x} , p \rangle_{\cH} > \langle x_i , p \rangle_{\cH}$ for all $i \in \llbracket 1,n \rrbracket$. Whence, it holds that 
\begin{equation*}
\big\langle \nabla_y \cJ_L^{\blambda}(\bx,\bar{x} \,|\,\ba) , p \big\rangle = \sum_{i=1}^n \underset{>0}{\underbrace{\lambda_i \, a_i \cdot \|\bar{x} - x_i\|_\cH^{a_i - 2}}} \; \langle \bar{x} - x_i , p \rangle > 0,
\end{equation*}
which contradicts the optimality condition $\nabla \cJ(\bar{x}) = 0$. The first claim made in \eqref{eq:convGMed} is a straightforward applications of~\eqref{eq:convHull} that we just proved, whereas for the second, we only need to check that $\QrGm_L^{\blambda}(\bx\,|\,x_0)$ belongs to the convex hull.
\end{proof}

\begin{rem}[Failure of representer theorems for medians in Banach spaces]
It is worth recalling that the conclusions of Lemma~\ref{lem:Hull} fail to hold if the ambient Hilbert space $\cH$ is replaced by a Banach space $\cX$, even in the ideal case in which the latter is, say, strictly convex and reflexive with strictly convex dual. Indeed, consider $\R^d$ with $d >2$ endowed with the $\ell^p$-norm
\begin{equation*}
\|x\|_p = \Big( |x_1|^p + \ldots + |x_d|^p \Big)^{\nicefrac{1}{p}}    
\end{equation*}
for some $p>1$, and let $\cJ(\boldsymbol{e},y) = \|y - e_1\|_p + \ldots + \|y - e_d\|_p$ where $\boldsymbol{e} \in (\R^d)^d$ is a canonical basis. Then $\bar{x} \in \argmin_{y \in \R^d} \cJ(\boldsymbol{e},y)$ is invariant under coordinate permutations, and there is $\bar{t} \in \R$ for which 
\begin{equation*}
\bar{x} = \bar{t} \cdot \sum_{i=1}^d e_i, 
\end{equation*}
whence
\begin{equation*}
\cJ(\bar{x}) = \sum_{i=1}^d \bigg\| \, e_i - \bar{t} \cdot \sum_{j=1}^d e_j \, \bigg\| = d \Big( |1-\bar{t}|^p + (d-1)|\bar{t}|^p \Big)^{\nicefrac{1}{p}}.
\end{equation*}
In particular, the point $\bar{t} \in \R$ must be a minimizer of the convex function $f(t) := |1-t|^p + (d-1)|t|^p$, and therefore satisfy
\begin{equation*}
f'(\bar{t}) = -p(1-\bar{t})^{p-1} + (d-1)p \, \bar{t}^{p-1} = 0
\end{equation*}
which by elementary algebraic manipulations can be shown to amount to 
%
%
\begin{equation*}
\left(\frac{1-\bar{t}}{\bar{t}}\right)^{p-1} = d-1.
\end{equation*}
Yet, the condition that $\bar{x} \in \conv(\bx)$ can only hold provided that $\bar{t}=\nicefrac{1}{d} \in (0,1)$. Since we picked $d>2$ to avoid trivial medians, this together with the above identity forces $p=2$, which corresponds to the Euclidean setting.
\end{rem}


\section{Quadratically regularized Wasserstein medians}
\label{section:QrWm}

In this section, we introduce meaningful quadratic regularizations of median-type aggregators in Wasserstein spaces, which fit the general resilience framework of Section \ref{section:Resilient} while being amenable to numerical simulations.

As highlighted in the statement and proof of Theorem \ref{thm:Resilience} above, the key ingredient for H\"older resilience is the validity of a uniform convexity estimate on the minimized cost, which in practice stems from the uniform convexity of the quadratic regularizer along perturbations for which the median terms are themselves convex. Our modeling choices were guided by the following convexity properties of Wasserstein distances, for which we refer to \cite[Chapter 7]{AGS}, \cite[Section 7.2]{OTAM} and \cite{Delalande2023} (see also \cite[Proposition 2.3]{CarlierD2024}). 

\begin{prop}[Miscellaneous convexity properties in Wasserstein spaces]
\label{prop:ConvexityWass}
The following hold. 
\begin{enumerate}
\item[$(a)$] For any $p \in [1,+\infty)$ and $\rho \in \cP_p(\R^d)$, the map $\mu \mapsto W_p^p(\rho,\mu)$ is \emph{flatly convex}, namely 
\begin{equation*}
W_p^p \Big( \rho , (1-t)\mu + t \nu \Big) \leq (1-t) \cdot W_p^p(\rho,\mu) + t \cdot W_p^p (\rho,\nu), 
\end{equation*}
for every $\mu,\nu \in \cP_p(\R^d)$, and \emph{strictly convex} if $p >1$ and $\rho \in \cP_p^r(\R^d)$. On the other hand, the map $\mu \mapsto W_p(\rho,\mu)$ is typically nonconvex for $p \neq 1$. 
\item[$(b)$] For any compact set $K \subset \R^d$ and every $\rho \in \cP^{\infty}_c(\R^d)$, the map $\mu \in \cP(K) \mapsto \tfrac{1}{2} W_2^2(\rho,\mu)$ is \emph{flatly uniformly convex}, in the sense that there is a constant $C_{\rho,K,d} > 0$ as well as a function $\omega : [0,1] \to [0,1]$ satisfying
\begin{equation*}
\lim_{t \to 0^+} \frac{\omega(t)}{t} = \lim_{t \to 1^-} \frac{\omega(t)}{1-t} = 1, 
\end{equation*}
such that
\begin{equation*}
\tfrac{1}{2} W_2^2\Big( \rho , (1-t)\mu + t\nu \Big) \leq (1-t) \cdot \tfrac{1}{2} W_2^2(\rho,\mu) + t \cdot \tfrac{1}{2} W_2^2(\rho,\nu) - C_{\rho,K,d} \, \omega(t) \cdot \tfrac{1}{2} W_2^6(\mu,\nu)
\end{equation*}
for every $\mu,\nu \in \cP(K)$.
\item[$(c)$] For every $\rho \in \cP_2^{\mathrm{ac}}(\R^d)$ and each $\sigma \in \cP_2(\R^d)$, the map $\mu \in \cP_2(\R^d) \mapsto W_{2,\rho}(\sigma,\mu)$ is \emph{convex along generalized geodesics with base $\rho$}, in the sense that 
\begin{equation*}
W_{2,\rho} \Big( \sigma \, , \, \big( (1-t) \cdot T_{\rho}^{\mu} + t \cdot T_{\rho}^{\nu} \big)_{\sharp} \rho \Big) \leq (1-t) \cdot W_{2,\rho}(\sigma,\mu) + t \cdot W_{2,\rho}(\sigma,\nu) 
\end{equation*}
for every $\mu,\nu \in \cP_2(\R^d)$.
\item[$(d)$] For every $\rho \in \cP_2^{\mathrm{ac}}(\R^d)$, the map $\mu \in \cP_2(\R^d) \mapsto \tfrac{1}{2} W_2^2(\rho,\mu)$ is \emph{strongly convex along generalized geodesics with base $\rho$}, namely 
\begin{equation*}
\tfrac{1}{2} W_2^2 \Big( \rho \, , \, \big( (1-t) \cdot T_{\rho}^{\mu} + t \cdot T_{\rho}^{\nu} \big)_{\sharp} \rho \Big) \leq (1-t) \cdot \tfrac{1}{2} W_2^2(\rho,\mu) + t \cdot \tfrac{1}{2} W_2^2(\rho,\nu) - t (1-t) \cdot \tfrac{1}{2} W_2^2(\mu,\nu) 
\end{equation*}
for every $\mu,\nu \in \cP_2(\R^d)$.
\end{enumerate}
\end{prop}

\begin{rem}[On the implications of Proposition \ref{prop:ConvexityWass}]
Points $(d)$ and $(b)$ above tell us that quadratic regularizers in Wasserstein spaces are strongly convex along suitably chosen geodesics starting from a proper base measure, and uniformly convex along flat interpolations provided the base measure is sufficiently regular and the target ones uniformly compactly supported. Complementarily, item $(a)$ says that the only Wasserstein metric that is flatly convex without being raised to any power is $W_1$, whereas item $(c)$ suggests that trading classical $W_2$-metrics for their linearizations around $\rho \in \cP_2^{\mathrm{ac}}(\R^d)$ opens the door to convexity along generalized geodesics. This observation is crucial, as the usual Wasserstein distance $\mu \in \cP_2(\R^d) \mapsto W_2(\mu,\sigma)$ is typically neither convex nor semiconvex along generalized geodesics, see e.g. \cite[Section 7.3]{AGS}.   
\end{rem}


\subsection{Quadratically regularized \texorpdfstring{$W_1$-medians}{}}
\label{subsection:QrWm}

In this section, we study quadratic regularization of geometric medians defined using classical Wasserstein spaces. These aggregators, which were introduced quite recently in \cite{Carlier2024,You2025} for $W_1$- and $W_2$-metrics respectively, are defined as follows.  

\begin{dfn}[Wasserstein medians]
Given a collection of weights $\blambda \in \Delta^n$, the \emph{$W_p$-median} of a family of measures $\bnu \in \cP_p(\R^d)^n$ is defined by 
\begin{equation}
\textsc{WMed}_p^{\blambda}(\bnu) = \argmin_{\mu \in \cP_p(\R^d)} \sum_{i=1}^n \lambda_i \, W_p(\nu_i,\mu).
\end{equation}
The latter aggregator is typically not single-valued, and as noted in \cite{Carlier2024}, the above problem is convex only in the particular case in which $p =1$, see also Theorem \ref{thm:Resilience}-$(a)$. 
\end{dfn}

\begin{rem}[Sharp breakdown point of $W_p$-medians]
As previously mentioned in Section \ref{section:Preli}, it was shown in \cite{Carlier2024} that 
\begin{equation*}
\bbeta^n_{W_p}(\WMed_p^{\blambda}) = \min \bigg\{ \sum_{j \in J} \lambda_j \; \Big| \; J \subset \llbracket 1,n \rrbracket \text{ and } \sum_{j\in J} \lambda_j \geq \frac{1}{2} \bigg\} \qquad \text{and} \qquad \bbeta^{\infty}_{W_p}(\WMed_p) = \frac{1}{2}
\end{equation*}
for each $p \in [1,+\infty)$. In Proposition \ref{prop:breakdown}, we proved a metric version of \cite[Theorem 2.1]{Lopuhaa1991}, which states that such a breakdown point is optimal in the class of aggregators over metric spaces $(\cX,\dist)$ which are equivariant under the action of a group of unbounded isometries $(\G,\circ)$. Such is the case for $(\cP_p(\R^d),W_p)$ for any $p \in [1,+\infty)$, which admits a large score of isometries (see e.g. the detailed treatment of the case $p=2$ in \cite{Kloeckner2010}), although for our purposes, it suffices to consider $(\R^d,+)$ whose induced actions $g : \cP_p(\R^d) \to \cP_p(\R^d)$ are defined by 
\begin{equation*}
g^h \cdot \mu := (\Id+h)_{\sharp} \mu.  
\end{equation*}
One may easily check that those are indeed group actions over $\cP_p(\R^d)$, since 
\begin{equation*}
(g^{h_1} \circ g^{h_2})_{\sharp} \mu = (g^{h_1+h_2})_{\sharp} \mu = g^{h_1}_{\sharp}(g^{h_2}_{\sharp} \mu).
\end{equation*}
The latter are also isometries, due to the fact that
\begin{equation*}
\Gamma(g^h \cdot  \mu,g^h \cdot  \nu) = \Big\{ (g^h \circ \pi^1,g^h \circ \pi^2)_{\sharp} \gamma \; \big|\; \gamma \in \Gamma(\mu,\nu) \Big\}. 
\end{equation*}
for every $\mu,\nu \in \cP_p(\R^d)$, which upon taking $\gamma^{\star} \in \Gamma_o(\mu,\nu)$ and $\gamma^{\star}_h \in \Gamma_o(g^h_{\sharp}\mu,g^h_{\sharp}\nu)$ entails that 
\begin{equation*}
\left\{
\begin{aligned}
& W_p(g^h \cdot  \mu,g^h \cdot  \nu)^p \leq \int_{\R^{2d}} \|x-y\|\, \mathrm{d} (g^h \circ \pi^1,g^h \circ \pi^2)_{\sharp} \gamma^{\star}(x,y) = W_p(\mu,\nu)^p \\
& W_p(\mu,\nu)^p \leq \int_{\R^{2d}} \|x-y\|\, \mathrm{d} (g^h \circ \pi^1,g^h \circ \pi^2)_{\sharp}^{-1} \gamma^{\star}_h(x,y) \leq W_p(g^h \cdot  \mu,g^h \cdot  \nu)^p
\end{aligned}
\right.
\end{equation*}
so that $W_p(g^h \cdot  \mu,g^h \cdot  \nu) = W_p(\mu,\nu)$. Such isometries are also clearly unbounded, since 
\begin{equation*}
W_p(\mu,g^h \cdot  \mu) \geq \cM_p(g^h \cdot \mu) - \cM_p(\mu) ~\underset{\|h\| \to +\infty}{\longrightarrow}~ +\infty
\end{equation*}
by the triangle inequality. 
\end{rem}

\begin{dfn}[Quadratically regularized $W_1$-median]
Given a collection of weights $\blambda \in \Delta^n$ and $L > 0$, the \emph{quadratically regularized $W_1$-median} with anchor $\nu_0 \in \cP_2(\R^d)$ of a family of measures $\bnu \in \cP_2(\R^d)^n$ is defined by 
\begin{equation}
\label{eq:QrWMed}
\QrWm^{\blambda}_{L}(\bnu \, |\, \nu_0) = \argmin\limits_{\mu \in \cP_2(\R^d) } \frac{1}{L} \cdot W_2^2(\nu_0,\mu) + \sum_{i=1}^n \lambda_i \cdot W_{1}(\nu_i,\mu).
\end{equation}
\end{dfn}

It is quite clear from Proposition \ref{prop:ConvexityWass}-$(a)$ above that $\QrWm^{\blambda}_{L}(\bnu \, |\, \nu_0)$ is a singleton as soon as $\nu_0 \in \cP_2^{\mathrm{ac}}(\R^d)$. Moreover, much like classical Wasserstein barycenters \cite{Agueh2011} and $W_1$-medians \cite{Carlier2024}, the above aggregator admits an equivalent definition in terms of a multi-marginal optimal transport problem. 

\begin{prop}[Multi-marginal reformulation of quadratically regularized $W_1$-medians]
\label{prop:MMOT}
Every element $\bar{\mu}_L \in \QrWm^{\blambda}_{L}(\bnu \, |\, \nu_0)$ is of the form $(T_L^{\blambda})_{\sharp} \bar{\bgamma}$, where 
\begin{subequations}
\begin{equation}
\label{eq:Tlambda}
T_L^{\blambda}(\bx) = \argmin_{y \in \R^d} \frac{1}{L} \cdot \|y-x_0\|^2 + \sum_{i=1}^n \lambda_i \cdot \|y - x_i\|
\end{equation}
and
\begin{equation}
\label{eq:MMOT}
\begin{aligned}
\bar{\bgamma} \in \argmin \Bigg\{ \int_{(\R^d)^{n+1}} \bigg( \frac{1}{L} \cdot \|T_L^{\blambda}(\bx)-x_0\|^2 + \sum_{i=1}^n \lambda_i \cdot \|T_L^{\blambda}(\bx) - x_i\|  \bigg) \, \mathrm{d} \bgamma(\bx) \; \bigg| \; \bgamma \in \Gamma(\nu_0,\ldots,\nu_n) \Bigg\}. 
\end{aligned}
\end{equation}
\end{subequations}
\end{prop}

\begin{proof}
The argument follows a standard blueprint setup in \cite{Agueh2011} and exploited in \cite{Carlier2024} for Wasserstein medians. First, given any $\bar{\mu}_L \in \QrWm_L^{\blambda}(\bnu\,|\,\nu_0)$, pick a family of optimal plans $\bar{\gamma}^i \in \Gamma_o(\bar{\mu}_L,\nu_i)$ for $i \in \llbracket 0,n \rrbracket$ and by the gluing lemma (see e.g. \cite[Lemma 5.3.2]{AGS}) define the plan $\bar{\bmu} \in \Gamma(\nu_0,\ldots,\nu_n,\bar{\mu}_L)$ as 
\begin{equation*}
\bar{\bmu}_L := \int_{\R^d} \, \big( \bar{\gamma}^0_y \otimes \ldots \otimes \bar{\gamma}^n_y \big)\, \mathrm{d} \bar{\mu}_L(y)
\end{equation*}
where $\{\bar{\gamma}^i_y\}_{y \in \R^d} \subset \cP(\R^d)$ are the disintegrations of the optimal plans against their shared marginal $\bar{\mu}_L$. Then, letting $\bar{\bgamma}_L := (\pi^0,\ldots,\pi^n)_{\sharp} \bar{\bmu}_L \in \Gamma(\nu_0,\ldots,\nu_n)$ it holds by construction that 
\begin{equation}
\label{eq:MMOTProof1}
\begin{aligned}
& \frac{1}{L} \cdot W_2^2(\bar{\mu}_L,\nu_0) + \sum_{i=1}^n \lambda_i \cdot W_{1}(\bar{\mu}_L,\nu_i) \\
& \hspace{0.45cm} = \int_{(\R^d)^{n+2}} \bigg( \frac{1}{L} \cdot \|y-x_0\|^2 + \sum_{i=1}^n \lambda_i \cdot \|y - x_i\|  \bigg) \, \mathrm{d} \bmu(\bx,y) \\
& \hspace{0.45cm} \geq \int_{(\R^d)^{n+1}} \bigg( \frac{1}{L} \cdot \|T_L^{\blambda}(\bx)-x_0\|^2 + \sum_{i=1}^n \lambda_i \cdot \|T_L^{\blambda}(\bx) - x_i\|  \bigg) \, \mathrm{d} \bar{\bgamma}(\bx) 
\end{aligned}
\end{equation}
where the inequality is obtained through pointwise minimization in the integral using \eqref{eq:Tlambda}. Conversely, if $\bar{\bgamma} \in \Gamma(\nu_0,\ldots,\nu_n)$ is a minimizer of \eqref{eq:MMOT} and $\bar{\mu}_L := (T_L^{\blambda})_{\sharp} \bar{\bgamma}$, then necessarily $\gamma^i := (\pi^i,T_L^{\blambda})_{\sharp} \bar{\bgamma} \in \Gamma(\nu_i,\bar{\mu}_L)$ for each $ i \in \llbracket 0,n \rrbracket$, so that 
\begin{equation}
\label{eq:MMOTProof2}
\begin{aligned}
& \int_{(\R^d)^{n+1}} \bigg( \frac{1}{L} \cdot \|T_L^{\blambda}(\bx)-x_0\|^2 + \sum_{i=1}^n \lambda_i \cdot \|T_L^{\blambda}(\bx) - x_i\|  \bigg) \, \mathrm{d} \bar{\bgamma}(\bx) \\
& \hspace{0.45cm} = \frac{1}{L} \cdot \int_{\R^{2d}}  \|y-x_0\|^2 \, \mathrm{d} \gamma^0(x_0,y) + \sum_{i=1}^n \lambda_i \cdot \int_{\R^{2d}} \|y - x_i\| \, \mathrm{d} \gamma^i(x_i,y) \\
& \hspace{0.45cm} \geq \frac{1}{L} \cdot W_2^2(\bar{\mu}_L,\nu_0) + \sum_{i=1}^n \lambda_i \cdot W_{1}(\bar{\mu}_L,\nu_i)
\end{aligned}
\end{equation}
Thus, we infer from combining \eqref{eq:MMOTProof1} and \eqref{eq:MMOTProof2} that both the minimization problems have the same value. From there, it is not hard to show by repeating the above constructions that $\bar{\mu}_L \in \QrWm_L^{\blambda}(\bnu \,|\,\nu_0)$ if and only if $\bar{\mu}_L = (T_L^{\blambda})_{\sharp} \bar{\bgamma}$ for some minimizer $\bar{\bgamma} \in \Gamma(\nu_0,\ldots,\nu_n)$ of \eqref{eq:MMOT}. 
\end{proof}

In addition to its relevance for numerics, 
the above multi-marginal representation allows to show that moment or support bounds on the anchor $\nu_0 \in \cP(\R^d)$ are automatically inherited by the aggregator. 

\begin{cor}[Uniform moment bounds on regularized medians]
\label{cor:MomentControl}
Suppose that $\nu_0 \in \cP_p(\R^d)$ for some $p \in [2,+\infty)$. Then necessarily $\QrWm_L^{\blambda}(\bnu\,|\,\nu_0) \in \cP_p(\R^d)$. Similarly if $\nu_0 \in \cP_c(\R^d)$, there exists a compact set $K_{\nu_0} \subset \R^d$ such that $\QrWm_L^{\blambda}(\bnu\,|\,\nu_0) \in \cP(K_{\nu_0})$. 
\end{cor}

\begin{proof}
As shown in Proposition \ref{prop:MMOT}, every $\bar{\mu}_L \in \QrWm_L^{\blambda}(\bnu \,|\, \nu_0)$ is of the form $\bar{\mu}_L = (T_L^{\blambda})_{\sharp}\bar{\bgamma}$ where 
\begin{equation*}
T_L^{\blambda}(\bx) = \argmin_{y \in \R^d} \frac{1}{L} \cdot \|y-x_0\|^2 + \sum_{i=1}^n \lambda_i \cdot \|y - x_i\|. 
\end{equation*}
The above problem being strongly convex, it follows from standard sensitivity results (see e.g. \cite[Section 4.4.1]{BonnansShapiro}) that $\bx \in (\R^d)^{n+1} \mapsto T_L^{\blambda}(\bx) \in \R^d$ is continuous. Then, by the basic Lagrange multiplier rule for convex optimization problems (see e.g. \cite[Theorem 8.15]{Rockafellar1998}, there is a collection of vectors $(v_1(\bx),\ldots,v_n(\bx)) \in (\R^d)^n$ satisfying $v_i(\bx) \in \partial (\|\cdot\|)(T_L^{\blambda}(\bx)-x_i)$, and such that 
\begin{equation*}
T_L^{\blambda}(\bx) - x_0 + \frac{L}{2} \cdot \sum_{i=1}^n \lambda_i v_i(\bx) = 0. 
\end{equation*}
Since $\|v_i(\bx)\| \leq 1$ for every $i \in \{1,\ldots,n\}$ and $\blambda \in \Delta^n$, this allows us to derive the moment estimate
\begin{equation*}
\begin{aligned}
\cM_p(\bar{\mu}_L) & = \bigg( \int_{(\R^d)^{n+1}} \|T_L^{\blambda}(\bx) \|^p \, \mathrm{d} \bar{\bgamma}(\bx) \bigg)^{\nicefrac{1}{p}} \\
& \leq \bigg( \int_{(\R^d)^{n+1}} \|x_0\|^p \, \mathrm{d} \bar{\bgamma}(\bx) \bigg)^{\nicefrac{1}{p}} + \frac{L}{2} \cdot \sum_{i=1}^n \lambda_i \cdot \bigg( \int_{(\R^d)^{n+1}} \|v_i(\bx)\|^p \, \mathrm{d} \bar{\bgamma}(\bx) \bigg)^{\nicefrac{1}{p}} \\
& \leq \cM_p(\nu_0) + \frac{L}{2}, 
\end{aligned}
\end{equation*}
which implies that $\bar{\mu}_L \in \cP_p(\R^d)$ as soon as $\nu_0 \in \cP_p(\R^d)$. It follows from similar considerations that 
\begin{equation*}
\| T_L^{\blambda}(\bx) \| \leq \|x_0\| + \frac{L}{2}
\end{equation*}
so that $\supp(\bar{\mu}_L) \subset K_{\nu_0} := (\supp(\nu_0) + B(0,L/2))$ is compactly supported when $\nu_0 \in \cP_c(\R^d)$.  
\end{proof}

In the following theorem, we show that the above aggregator is H\"older resilient in the sense of Definition \ref{def:Resilience}.

\begin{thm}[Hölder resilience of quadratically regularized $W_1$-medians]
Suppose that $\nu_0 \in \cP_c^{\infty}(\R^d)$. Then for every $\blambda \in \Delta^n$, the aggregator 
\begin{equation*}
\QrWm^{\blambda}_L(\bullet \,|\, \nu_0) : \cP_2(\R^d)^n \to \cP_c(\R^d)    
\end{equation*}
is H\"older resilient with exponent $\alpha = \nicefrac{1}{5}$.  
\end{thm}

\begin{proof}
The above result follows from an application of the general criterion for H\"older resilience derived in Theorem \ref{thm:Resilience}. To lighten the notations given $\bnu \in \cP_2(\R^d)^n$, we let $\bar{\mu}_L := \QrWm^{\blambda}_L(\bnu \,|\, \nu_0)$, which belongs to $\cP(K_{\nu_0})$ with $K_{\nu_0} := (\supp(\nu_0) + B(0,L/2))$ by Corollary \ref{cor:MomentControl}. Recalling the definition \eqref{eq:QrWMed} of the aggregator, we are in the framework of the aforementioned proposition with 
\begin{equation*}
\cJ_L^{\blambda}(\bnu,\mu) := \frac{1}{L} \cdot W_2^2(\nu_0,\mu) + \sum_{i=1}^n \lambda_i \cdot W_{1}(\nu_i,\mu),
\end{equation*}
and note that point $(i)$ thereof is clearly satisfied with $(\cX,\dist) = (\cP_2(K),W_2)$ and $\dsf_i = W_1$ for each $i \in \llbracket 1,n \rrbracket$. Concerning point $(ii)$, it follows from items $(a)$ and $(b)$ of Proposition \ref{prop:ConvexityWass} above combined with the fact that $\bar{\mu}_L = \argmin_{\mu \in \cP_2(\R^d)} \cJ_L^{\blambda}(\bnu,\mu)$ that 
\begin{equation*}
\begin{aligned}
\cJ_L^{\blambda}(\bnu,\bar{\mu}_L) & \leq \cJ_L^{\blambda}\Big( \bnu , (1-t)\bar{\mu}_L + t \mu \Big) \\
& = \frac{1}{L} \cdot \underset{\text{Proposition \ref{prop:ConvexityWass}-$(b)$}}{\underbrace{W_2^2 \Big( \nu_0 \,, (1-t) \bar{\mu}_L + t \mu  \Big)}} + \sum_{i=1}^n \lambda_i \cdot \underset{\text{Proposition \ref{prop:ConvexityWass}-$(a)$}}{\underbrace{W_1 \Big( \nu_i \,, (1-t) \bar{\mu}_L + t \mu \Big)}} \\
& \leq (1-t) \cdot \cJ_L^{\blambda}(\bnu,\bar{\mu}_L) + t \cdot \cJ_L^{\blambda}(\bnu,\mu) - \frac{C_{\nu_0,d}}{L} \, \omega(t) \cdot W_2^6(\bar{\mu}_L,\mu)
\end{aligned}
\end{equation*}
for all $t \in [0,1]$ and every $\mu \in \cP(K_{\nu_0})$, where $C_{\nu_0,d} := C_{\nu_0,K_{\nu_0},d}$ given as in Proposition \ref{prop:ConvexityWass}-$(b)$. Upon rearranging terms, dividing by $t > 0$ and finally letting $t \to 0^+$, we finally infer that 
\begin{equation*}
\cJ_L^{\blambda}(\bnu,\mu) \geq \cJ_L^{\blambda}(\bnu,\bar{\mu}_L) + \frac{C_{\nu_0,d}}{L} \cdot W_2^6(\bar{\mu}_L,\mu)
\end{equation*}
for every $\mu \in \cP(K_{\nu_0})$. There lastly remains to observe that to due to Corollary \ref{cor:MomentControl}, any candidate minimizer associated with another collection of inputs $\tilde{\bnu} \in \cP_2(\R^d)^n$ is supported in $K_{\nu_0}$ to conclude using following Theorem \ref{thm:Resilience} that $\QrWm_L^{\blambda}(\bullet \,|\,\nu_0) : \cP_2(\R^d)^n \to \cP_c(\R^d)$ is $\nicefrac{1}{5}$-H\"older resilient.
\end{proof}


\subsection{Quadratically regularized \texorpdfstring{$W_{2,\nu_0}$-medians}{}}
\label{subsection:QrLWm}

In this section, we shift our focus and consider quadratic regularizations of a new median-type aggregator defined in terms of linearized Wasserstein metrics $W_{2,\nu_0}$ around an anchor $\nu_0 \in \cP_2^{\mathrm{ac}}(\R^d)$.  

\begin{dfn}[Linearized Wasserstein medians]
Given a collection of weights $\blambda \in \Delta^n$ and an anchor $\nu_0 \in \cP_2^{\mathrm{ac}}(\R^d)$, the \emph{$W_{2,\nu_0}$-median} of a family of measures $\bnu \in \cP_2(\R^d)^n$ is defined by 
\begin{equation}
\label{eq:LinearizedMedian}
\LWMed^{\blambda}_2(\bnu\,|\,\nu_0) = \argmin_{\mu \in \cP_2(\R^d)} \sum_{i=1}^n \lambda_i \cdot W_{2,\nu_0}(\nu_i,\mu).
\end{equation}
\end{dfn}

Owing to the consequence of Brenier's theorem recollected in Theorem \eqref{thm:Brenier}, it can be shown that linearized Wasserstein medians can be expressed as pushforwards of $\nu_0 \in \cP_2^{\mathrm{ac}}(\R^d)$ through the (Hilbertian) median of the family $\bT_{\nu_0}^{\bnu} := (T_{\nu_0}^{\nu_1},\ldots,T_{\nu_0}^{\nu_n}) \in L^2(\R^d,\R^d;\nu_0)^n$ of input transport maps.

\begin{prop}[Structure and representation of $W_{2,\nu_0}$-medians]
\label{prop:LinMed}
Given a collection of weights $\blambda \in \Delta^n$, one has that $\bar{\mu} \in \LWMed^{\blambda}_2(\bnu\,|\,\nu_0)$ \emph{if and only if} it is of the form $\bar{\mu} = (\bar{T}_{\nu_0})_{\sharp} \nu_0$, where 
\begin{equation}
\label{eq:BrenierMedian}
\bar{T}_{\nu_0} \in \GMed^{\blambda}(\bT_{\nu_0}^{\bnu}) = \argmin \bigg\{ \sum_{i=1}^n \lambda_i \, \|T-T_{\nu_0}^{\nu_i} \|_{L^2(\nu_0)} \; \Big| \; T \in L^2(\R^d,\R^d;\nu_0) \bigg\}.     
\end{equation}
In particular, $\LWMed^{\blambda}_2(\bnu\,|\,\nu_0)$ is a singleton as soon as the elements of $\bT_{\nu_0}^{\bnu}$ are not all collinear. Moreover, the maps $\bar{T}_{\nu_0}$ are of the form
\begin{equation}
\label{eq:ConvexHullbarT}
\bar{T}_{\nu_0} = \sum_{i=1}^n \bar{\theta}_i(\bT_{\nu_0}^{\bnu}) \cdot T_{\nu_0}^{\nu_i},
\end{equation}
where the weights $\bar{\btheta}(\bT_{\nu_0}^{\bnu}) \in \Delta^n$ are characterized through the minimization problem
\begin{equation}
\label{eq:WeightsbarT}
\bar{\btheta}(\bT_{\nu_0}^{\bnu}) \in \argmin_{\btheta \in \Delta^n} \sum_{i=1}^n \lambda_i \cdot \, \bigg\|\, \sum_{j=1}^n \theta_j \cdot T_{\nu_0}^{\nu_j} -T_{\nu_0}^{\nu_i} \,\bigg\|_{L^2(\nu_0)}.
\end{equation}
Finally, if $\nu_i \in \cP_2^{\mathrm{ac}}(\R^d)$ for every $i \in \llbracket 1,n \rrbracket$, then $\LWMed_2^{\blambda}(\bnu) \subset \cP_2^{\mathrm{ac}}(\R^d)$ as well.  
\end{prop}

\begin{proof}
From Theorem \ref{thm:Brenier}, every candidate measure $\mu \in \cP_2(\R^d)$ in \eqref{eq:LinearizedMedian} may be written as $\mu = (T_{\nu_0}^{\mu})_{\sharp} \nu_0$ for a uniquely defined map $T_{\nu_0}^{\mu} \in L^2(\R^d,\R^d;\nu_0) \cap \nabla \cC(\R^d)$, whence
\begin{equation*}
\min_{\mu \in \cP_2(\R^d)} \sum_{i=1}^n \lambda_i \cdot W_{2,\nu_0}(\nu_i,\mu) = \min \bigg\{ \sum_{i=1}^n \lambda_i \cdot \|T-T_{i}\|_{L^2(\nu_0)} \; \Big| \; T \in L^2(\R^d,\R^d;\nu_0) \cap \nabla \cC(\R^d) \bigg\}
\end{equation*}
and $\bar{\mu} \in \LWMed_2^{\blambda}(\bnu\,|\,\nu_0)$ if and only if $\bar{\mu} = (\bar{T}_{\nu_0})_{\sharp} \nu_0$ for a minimizer $\bar{T}_{\nu_0}$ of the above problem. Yet, since $\bT_{\nu_0}^{\bnu} \in (L^2(\R^d,\R^d;\nu_0) \cap \nabla\cC(\R^d))^n$ by definition, it follows from Lemma \ref{lem:Hull} that 
\begin{equation}
\label{eq:ConvexHullProof}
\GMed^{\blambda}(\bT) \subset \conv(\bT) 
\end{equation}
which, upon recalling that $\nabla \cC(\R^d)$ is convex, implies that 
\begin{equation*}
\GMed^{\blambda}(\bT_{\nu_0}^{\bnu}) = \argmin \bigg\{ \sum_{i=1}^n \lambda_i \cdot \|T-T_{i}\|_{L^2(\nu_0)} \; \Big| \; T \in L^2(\R^d,\R^d;\nu_0) \cap \nabla \cC(\R^d) \bigg\}. 
\end{equation*}
Then, the uniqueness of the linearized Wasserstein median follows then from that of $\GMed^{\blambda}(\bT_{\nu_0}^{\bnu})$, which is itself guaranteed by the non-collinearity assumption, see e.g. \cite{Minsker2015} or the original proof in \cite{Kemperman1987}. To close the proof, we observe that the identity \eqref{eq:ConvexHullbarT} and the characterization \eqref{eq:WeightsbarT} of the weights follow straightforwardly from \eqref{eq:ConvexHullProof}, and that the absolute continuity statement may be proven exactly as \cite[Theorem 3.1]{Alvarez2016}. 
\end{proof}

\begin{rem}[Breakdown point of the $W_{2,\nu_0}$-median]
As we just showed that the $W_{2,\nu_0}$-median is actually the Hilbert median of the input Brenier maps whenever $\nu_0 \in \cP_2^{\mathrm{ac}}(\R^d)$, it directly follows   
\begin{equation*}
\bbeta_{W_{2,\nu_0}}^n(\LWMed^{\blambda}(\bullet\,|\,\nu_0)) = \min \bigg\{ \sum_{j \in J} \lambda_j \; \Big|\; J \subset \llbracket 1,n \rrbracket ~~\text{and}~~ \sum_{j \in J} \lambda_j \geq \frac{1}{2} \bigg\}
\end{equation*}
for every $\blambda \in \Delta^n$, and likewise $\bbeta_{W_{2,\nu_0}}^{\infty}(\LWMed^{\blambda}(\bullet\,|\,\nu_0)) = \nicefrac{1}{2}$. Moreover, $W_{2,\nu_0}$-medians incur the same relative sensitivity to adversarial inputs as Euclidean ones, including the lack of strategy-proofness exhibited in \cite{ElMhamdi2023}. 
\end{rem}

\begin{dfn}[Quadratically regularized $W_{2,\nu_0}$-median]
Given a collection of weights $\blambda \in \Delta^n$ and $L > 0$, the \emph{quadratically regularized $W_{2,\nu_0}$-median} with anchor $\nu_0 \in \cP_2^{\mathrm{ac}}(\R^d)$ of a family of measures $\bnu \in \cP_2(\R^d)^n$ is defined by 
\begin{equation}
\QrLWm^{\blambda}_{L}(\bnu|\nu_0) := \underset{\mu \in \cP_2(\R^d)}{\mathrm{argmin}} \, \frac{1}{L} \cdot W_2^2(\nu_0,\mu) + \sum_{i=1}^n \lambda_i \cdot W_{2,\nu_0}(\nu_i,\mu).
\end{equation}
\end{dfn}

In the following theorem, we precisely characterize the quadratically regularized linearized Wassertein median in terms of the optimal transport maps $\bT_{\nu_0}^{\bnu} \in L^2(\R^d,\R^d;\nu_0)$ from the prior to the input data. 

\begin{thm}[Structure of quadratically regularized $W_{2,\nu_0}$-medians] 
\label{thm:QrLWm}
Suppose that $\nu_0 \in \cP_2^{\mathrm{ac}}(\R^d)$. Then for every $\blambda \in \Delta^n$, it holds that 
\begin{equation}
\label{eq:QrLWmQrGM}
\QrLWm^{\blambda}_{L}(\bnu \,|\,\nu_0) = \QrGm^{\blambda}_{L}(\bT_{\nu_0}^{\bnu} |\,\Id)_\sharp\nu_0.
\end{equation}
where
\begin{equation}
\label{eq:QrGMlambda}
\QrGm^{\blambda}_{L}(\bT_{\nu_0}^{\bnu} |\,\Id) = \argmin \Bigg\{ \frac{1}{L} \cdot \|T - \Id\|_{L^2(\nu_0)}^2 + \sum_{i=1}^n \lambda_i \cdot \|T - T_{\nu_0}^{\nu_i} \|_{L^2(\nu_0)} \; \bigg| \; T \in L^2(\R^d,\R^d;\nu_0) \Bigg\}.
\end{equation}
Moreover, the above aggregator is of the form 
\begin{equation}
\label{eq:WeightQrGm}
\QrGm^{\blambda}_{L}(\bT_{\nu_0}^{\bnu} |\,\Id) = \sum_{i=0}^n \theta_i^L(\bT_{\nu_0}^{\bnu}) \cdot T_{\nu_0}^{\nu_i},   
\end{equation}
where the weights $\bar{\btheta}^L(\bT_{\nu_0}^{\bnu}) \in \Delta^{n+1}$ are characterized through the minimization problem
\begin{equation}
\label{eq:FiniteQrGm}
\bar{\btheta}^L(\bT_{\nu_0}^{\bnu}) \in \argmin_{\btheta \in \Delta^{n+1}} \Bigg\{ \frac{1}{L} \cdot \bigg\| \sum_{i=0}^n \theta_i \cdot T_{\nu_0}^{\nu_i} - \Id \, \bigg\|^2_{L^2(\nu_0)} + \, \sum_{i=1}^n \lambda_i \cdot \, \bigg\|\, \sum_{j=0}^n \theta_j \cdot T_{\nu_0}^{\nu_j} -T_{\nu_0}^{\nu_i} \,\bigg\|_{L^2(\nu_0)} \Bigg\}. 
\end{equation} 
In particular if $\nu_0 \in \cP_p^{\mathrm{ac}}(\R^d)$ and $\bnu \in \cP_p^{\mathrm{ac}}(\R^d)^n$ for some $p \in [2,+\infty)$, then $\QrLWm_L^{\blambda}(\bnu\,|\,\nu_0) \in \cP_p^{\mathrm{ac}}(\R^d)$ as well, with 
\begin{equation}
\label{eq:MomentQrGm}
\cM_p \Big( \QrLWm_L^{\blambda}(\bnu\,|\,\nu_0) \Big) \leq \max\limits_{i \in \llbracket 0,n \rrbracket} \cM_p(\nu_i).
\end{equation}
\end{thm}

\begin{proof}
The demonstrations of \eqref{eq:QrLWmQrGM} through \eqref{eq:FiniteQrGm} follow the same path as those in Proposition \ref{prop:LinMed}, and similarly rely on the consequences of Brenier's theorem listed in Theorem \ref{thm:Brenier} and the convex hull lemma \ref{lem:Hull}. Concerning the moment bound \eqref{eq:MomentQrGm}, it suffices to observe that 
\begin{equation*}
\begin{aligned}
\cM_p \Big( \QrLWm_L^{\blambda}(\bnu\,|\,\nu_0) \Big) & = \bigg( \int_{\R^d} \|\bar{T}^L_{\nu_0}\|^p \, \mathrm{d} \nu_0(x) \bigg)^{\nicefrac{1}{p}} \\
& \leq \sum_{i=0}^n \theta_i^L(\bT_{\nu_0}^{\bnu}) \cM_p(\nu_i) \\
& \leq \max_{i \in \llbracket 0,n \rrbracket} \cM_p(\nu_i)
\end{aligned}
\end{equation*}
since $\btheta^L(\bT_{\nu_0}^{\bnu}) \in \Delta^{n+1}$. Finally, the fact $\QrLWm_L^{\blambda}(\bnu\,|\,\nu_0) \ll \Leb$ when the anchor and input measures are all absolutely continuous can be proven by combining \eqref{eq:QrLWmQrGM} and \eqref{eq:WeightQrGm} and repeating the arguments from \cite[Theorem 3.1]{Alvarez2016}, which themselves hinge upon the formula for the density of the pushforward, see e.g. \cite[Lemma 5.4.3]{AGS}.
\end{proof}

As a consequence of Theorem \ref{thm:QrLWm}, we infer the Lipschitz resilience of the above aggregator. 

\begin{cor}[Lipschitz resilience of quadratically regularized $W_{2,\nu_0}$-medians]
Suppose that $\nu_0 \in \cP_2^{\mathrm{ac}}(\R^d)$. Then for every $\blambda \in \Delta^n$, the aggregator $\QrLWm_L^{\blambda}(\bullet \,|\,\nu_0) : \cP_2(\R^d)^n \to \cP_2(\R^d)$ is Lipschitz resilient. 
\end{cor}

\begin{proof}
From \ref{thm:LipResGrGm}, we know that the aggregator $\QrGm^{\blambda}_{L}(\bullet \,|\,\Id) : L^2(\R^d,\R^d;\nu_0) \to L^2(\R^d,\R^d;\nu_0)$ is Lipschitz resilient with constant $L > 0$. Now since $\QrGm^{\blambda}_{L}(\bT_{\nu_0}^{\bnu} |\,\Id)$ is an optimal transport map for any $\bnu \in \cP_2(\R^d)^n$, this implies that 
\begin{equation*}
\begin{aligned}
W_2 \Big( \QrLWm^{\blambda}_{L}(\bnu\,|\,\nu_0), \QrLWm^{\blambda}_{L}(\sigma\,|\,\nu_0) \Big) & \leq W_{2,\nu_0} \Big( \QrLWm^{\blambda}_{L}(\bnu\,|\,\nu_0), \QrLWm^{\blambda}_{L}(\sigma\,|\,\nu_0) \Big) \\
& = \Big\| \, \QrGm^{\blambda}_{L}(\bT_{\nu_0}^{\bnu} |\,\Id) - \QrGm^{\blambda}_{L}(\bT_{\nu_0}^{\bsigma} |\,\Id) \, \Big\|_{L^2(\nu_0)} \\
& \leq L \cdot \CmD( \bT_{\nu_0}^{\bnu},\bT_{\nu_0}^{\bsigma}) \\
& = L \cdot \CmD(\bnu,\bsigma) 
\end{aligned}
\end{equation*}
for every $\bnu,\bsigma \in \cP_2(\R^d)^n$, whence the Lipschitz resilience. 
\end{proof}

Another interesting consequence of Theorem \ref{thm:QrLWm} is that the quadratically regularized linearized medians enjoy various densities bounds, depending on the absolutely continuity and moments of the anchor and input data. 

\begin{prop}[Density bounds for regularized $W_{2,\nu_0}$-medians]
\label{prop:Density}
The following estimates hold.
\begin{enumerate}
\item[$(a)$] Suppose that $\nu_0 = \rho_0 \cdot \Leb$ and $\nu_i = \rho_i \cdot \Leb$ for some $(\rho_0,\ldots,\rho_n) \in L^{\infty}(\R^d,\R_+)^{n+1}$. Then for every $L > 0$ and each $\blambda \in \Delta^n$, it holds that 
\begin{equation}
\label{eq:QrLWmDensity1}
\QrLWm_L^{\blambda}(\bnu\,|\,\nu_0) = \bar{\rho}_L(\bnu) \cdot \Leb
\end{equation}
where $\bar{\rho}_L(\bnu) \in L^{\infty}(\R^d,\R_+)$ satisfies $\|\bar{\rho}_L(\bnu)\|_{L^{\infty}} \leq \max\limits_{i \in \llbracket 0,n \rrbracket} \|\rho_i\|_{L^{\infty}}$. 
\item[$(b)$] Suppose that $\nu_0 = \rho_0 \cdot \Leb$ for some $\rho_0 \in L^{\infty}(\R^d,\R_+)$. Then for each $\blambda \in \Delta^n$ and every $\bnu \in \cP_2(\R^d)^n$, there exists $\bar{L} > 0$ such that for every $L \in (0,\bar{L})$, one has 
\begin{equation}
\label{eq:QrLWmDensity2}
\QrLWm_L^{\blambda}(\bnu\,|\,\nu_0) = \bar{\rho}_L^0 \cdot \Leb
\end{equation}
where $\bar{\rho}_L^0 \in L^{\infty}(\R^d,\R_+)$.  
\end{enumerate} 

\end{prop}

\begin{proof}
Both items are based on the representation \eqref{eq:QrLWmQrGM}, which ensures that $\QrLWm_L^{\blambda}(\bnu\,|\,\nu_0) =  (\bar{T}_{\nu_0}^L)_{\sharp} \nu_0$ where the map $\bar{T}_{\nu_0}^L \in L^2(\R^d,\R^d;\nu_0)$ is of the form 
\begin{equation}
\label{eq:BarThullProof}
\QrGm^{\blambda}_{L}(\bT_{\nu_0}^{\bnu} |\,\Id) = \sum_{i=0}^n \theta_i^L(\bT_{\nu_0}^{\bnu}) \cdot T_{\nu_0}^{\nu_i}
\end{equation}
for some weights $\bar{\btheta}^L(\bT_{\nu_0}^{\bnu}) \in \Delta^{n+1}$. Recalling that the inputs $\bT_{\nu_0}^{\bnu} \in L^2(\R^d,\R^d;\nu_0)^n$ are gradients of convex functions by Theorem \ref{thm:Brenier}, they are also Fréchet differentiable $\nu_0$-almost everywhere as a consequence of Aleksandrov's theorem (see e.g. \cite[Theorem 5.5.4]{AGS}). Then, because the input measures $\nu_i = \rho_i \cdot \Leb \in \cP_2^{\mathrm{ac}}(\R^d)$ are absolutely continuous with a bounded density for each $i \in \llbracket 1,n \rrbracket$, it stems from the necessary part of \cite[Lemma 5.5.3]{AGS} that
\begin{equation*}
\mathrm{det}(\nabla T_{\nu_0}^{\nu_i}(x)) > 0
\end{equation*}
for $\nu_0$-almost every $x \in \R^d$. Following the same arguments as in the proof of \cite[Theorem 3.1]{Alvarez2016}, it can be shown that $\bar{T}_{\nu_0}^L : \R^d \to \R^d$ is injective, which allows to apply the pushforward density formula of \cite[Lemma 5.5.3]{AGS} to obtain that \eqref{eq:QrLWmDensity1} holds with 
\begin{equation}
\label{eq:DensityQrLWm}
\bar{\rho}_L(\bnu) := \frac{\rho_0}{\mathrm{det}(\nabla \bar{T}_{\nu_0}^L)} \circ (\bar{T}_{\nu_0}^L)^{-1}. 
\end{equation}
To conclude, we combine the concavity of the $\mathrm{det}^{\nicefrac{1}{d}}(\bullet)$ over positive semidefinite matrices (see e.g. \cite[Lemma 5.23]{Villani1}) and the convexity of $t \in \R_+^{\star} \mapsto 1/t^d\in \R_+$ for $d \in \N$ to infer that
\begin{equation*}
\begin{aligned}
\frac{\rho_0}{\mathrm{det}(\nabla \bar{T}_{\nu_0}^L)} \circ (\bar{T}_{\nu_0}^L)^{-1} & = \frac{\rho_0}{\Big( \mathrm{det}\big(\sum_{i=0}^n \theta^L_i(\bT_{\nu_0}^{\bnu}) \cdot \nabla T_{\nu_0}^{\nu_i} \big)^{\nicefrac{1}{d}} \Big)^d} \circ (\bar{T}_{\nu_0}^L)^{-1} \\
& \leq \frac{\rho_0}{\Big( \sum_{i=0}^n \theta^L_i(\bT_{\nu_0}^{\bnu}) \cdot \mathrm{det}(\nabla T_{\nu_0}^{\nu_i})^{\nicefrac{1}{d}} \Big)^d} \circ (\bar{T}_{\nu_0}^L)^{-1} \\
& \leq \sum_{i=0}^n \theta^L_i(\bT_{\nu_0}^{\bnu}) \frac{\rho_0}{\mathrm{det}(\nabla T_{\nu_0}^{\nu_i})} \circ (\bar{T}_{\nu_0}^L)^{-1}. 
\end{aligned}
\end{equation*}
It is then enough to realise using again \cite[Lemma 5.5.3]{AGS} that 
\begin{equation*}
\Big\|\, \frac{\rho_0}{\mathrm{det}(\nabla T_{\nu_0}^{\nu_i})} \circ (\bar{T}_{\nu_0}^L)^{-1} \, \Big\|_{L^{\infty}} = \|\rho_i \circ T_{\nu_0}^{\nu_i} \circ (\bar{T}_{\nu_0}^L)^{-1} \|_{L^{\infty}} \leq \|\rho_i\|_{L^{\infty}}  
\end{equation*}
for each $i \in \llbracket 1,n \rrbracket$ to conclude since $\btheta^L(\bT_{\nu_0}^{\bnu}) \in \Delta^{n+1}$.

To prove item $(b)$, we follow the same the arguments but this time aim at showing that there exists $\bar{\btheta}^L(\bT_{\nu_0}^{\bnu}) \in \Delta^{n+1}$ realizing \eqref{eq:BarThullProof} and such that $\theta_0^L(\bT_{\nu_0}^{\bnu}) > 0$ for $L > 0$ small enough. If $\bar{T}_{\nu_0}^L = \Id$, we simply take $\theta_0^L(\bT_{\nu_0}^{\bnu})=1$. Otherwise, the only way one may not find $\theta^L_0(\bT_{\nu_0}^{\bnu}) > 0$ in this expression corresponds to the situation in which $\bar{T}_{\nu_0}^L = T_{\nu_0}^{\nu_j}$ for some $j \in \llbracket 1,n \rrbracket$. To exclude this possibility, we note that by writing the first-order optimality condition for $\QrGm^{\blambda}_{L}(\bT_{\nu_0}^{\bnu} |\,\Id)$, it can be shown through similar arguments as those discussed in Corollary \ref{cor:MomentControl} above that 
\begin{equation*}
\|\bar{T}_{\nu_0}^L - \Id \|_{L^2(\nu_0)} \leq \frac{L}{2}.
\end{equation*}
At this stage, if follows if we choose $\bar{L} > 0$ in such a way that 
\begin{equation*}
\bar{L} <\frac{1}{4} \min_{\substack{i \in \llbracket 1,n\rrbracket \\ T_{\nu_0}^{\nu_i} \neq \Id}} \|T_{\nu_0}^{\nu_i} - \Id \|_{L^2(\nu_0)}, 
\end{equation*}
then necessarily $\bar{T}_{\nu_0}^L \neq T_{\nu_0}^{\nu_i}$ for each $i \in \llbracket 1,n \rrbracket$ whenever $L \in (0,\bar{L})$, and the first-order condition 
\begin{equation*}
\bar{T}_{\nu_0}^L = \frac{\nicefrac{2}{L}\cdot \Id + \sum\limits_{i=1}^n \nicefrac{\lambda_i}{\|\bar{T}_{\nu_0}^L - T_{\nu_0}^{\nu_i} \|_{L^2(\nu_0)}} \cdot T_{\nu_0}^{\nu_i}}{\nicefrac{2}{L} +\sum\limits_{i=1}^n \nicefrac{\lambda_i}{\|\bar{T}_{\nu_0}^L - T_{\nu_0}^{\nu_i} \|_{L^2(\nu_0)}}}
\end{equation*}
provides an admissible representation with $\theta^L_0(\bT_{\nu_0}^{\bnu}) > 0$. Then, by repeating the above arguments, there exists a strongly convex function $\bar{\varphi}_{\nu_0}^L : \R^d \to \R \cup\{+\infty\}$ such that $\bar{T}_{\nu_0}^L = \nabla \bar{\varphi}_{\nu_0}^L$, whence the latter is injective and Fréchet differentiable almost everywhere, with
\begin{equation*}
\mathrm{det}(\nabla \bar{T}_{\nu_0}^L(x)) = \mathrm{det} \bigg( \theta^L_0(\bT_{\nu_0}^{\bnu}) \cdot \Id + \sum_{i=1}^n \theta^L_i(\bT_{\nu_0}^{\bnu}) \cdot \nabla T_{\nu_0}^{\nu_i}(x) \bigg) \geq \theta^L_0(\bT_{\nu_0}^{\bnu})^d
\end{equation*}
for $\nu_0$-almost every $x \in \R^d$. It then follows again from the pushforward density formula of \cite[Lemma 5.5.3]{AGS} that \eqref{eq:QrLWmDensity2} holds with 
\begin{equation*}
\bar{\rho}_L^0 := \frac{\rho_0}{\mathrm{det}(\nabla \bar{T}_{\nu_0}^L)} \circ (\bar{T}_{\nu_0}^L)^{-1}. 
\end{equation*}
In particular, one has that $\|\bar{\rho}_L^0\|_{L^{\infty}} \leq \|\rho_0\|_{L^{\infty}}/\theta^L_0(\bT_{\nu_0}^{\bnu})^d$, which concludes the proof. 
\end{proof}

\section{Majorization-minimization approximation of \texorpdfstring{$W_2$}{}-medians}
\label{subsection:MajorizeMinimize}

In this section, we outline a Majorization-Minimization (MM) scheme for the approximating stationary points of the $W_2$-median problem
\begin{equation}
\label{eq:MMW2Med}
\WMed_2^{\blambda}(\bnu) = \argmin_{\mu \in \cP_2(\R^d)} \sum_{i=1}^n \lambda_i \cdot W_2(\nu_i,\mu). 
\end{equation}
The latter relies on the properties of linearized Wasserstein distances, and on the density bounds derived above for quadratically regularized linearized medians. 

The term majorization-minimization broadly refers to a range of methods in which one approximaties either local minimima or critical points of a (not necessarily convex) Fréchet differentiable function $\cJ : \cX \to \R_+$ by inductively constructing a sequence $(x_k)_{k \in \N} \subset \cX$ through the update rule 
\begin{equation}
\label{eq:cost_MM}
x_{k+1} = \argmin_{y \in \cX} \, \cG(y\,|\,x_k). 
\end{equation}
Therein, the maps $\cG : \cX \times \cX \to \R_+$ are called \textit{surrogates}, and are usually taken as being jointly continuous and satisfying the following properties.   

\begin{dfn}[Admissible surrogates]
\label{def:Surrogates}
We say that $\cG(\bullet \,|\,x_k) : \cX \to \R_+$ is an \emph{admissible surrogate} for $\cJ$ at $x_k \in \cX$ if it satisfies the following properties. 
\begin{enumerate}
\item[$(i)$] It holds that $\cG(x_k\,|\,x_k) = \cJ(x_k)$ and $\cJ(y) \leq \cG(y\,|\,x_k)$ for every $y \in \cX$.
\item[$(ii)$] The maps $\cG(\bullet\,|\,x_k)$ comply with the strong convexity estimate
\begin{equation*}
\cG_L^{\blambda}(y\,|\,x_k) \geq \cG_L^{\blambda}(x_{k+1}\,|\,x_k) + \frac{1}{2L} \cdot \dist(x_{k+1},y)^2 
\end{equation*}
for every $y \in \cX$ and some $\nicefrac{1}{L}$-strongly convex 
for some constant $L >0$.
\end{enumerate}
\end{dfn}

Item $(i)$ in the above definition ensures that the sequence $(\cJ(x_k))_{k \in \N}$ is monotonically decreasing, whereas item $(ii)$ implies that the surrogates have a unique minimum $x_{k+1}$ which depend continuously on $x_k$. 
In our case, we let $(\cX,\dist) := (\cP_2(\R^d),W_2)$ and consider the nonconvex objective 
\begin{equation*}
\cJ^{\blambda}(\mu) := \sum_{i=1}^n \lambda_i \cdot W_2(\nu_i,\mu) 
\end{equation*}
defined for any given collections of weights $\blambda \in \Delta^n$ and input measures $\bnu \in \cP_2(\R^d)^n$. Given an anchor $\rho \in \cP_2^{\mathrm{ac}}(\R^d)$ and some $L >0$, we then introduce the maps 
\begin{equation*}
\cG^{\blambda}_L(\mu \,|\,\bnu,\rho) := \frac{1}{2L} \cdot W_2^2(\rho,\mu) + \sum_{i=1}^n \lambda_i \cdot W_{2,\rho}(\nu_i,\mu).
\end{equation*}
Owing to the basic estimate \eqref{eq:LinearizedWassIneq} between classical and linearized Wasserstein metrics, combined with the convexity estimates of Proposition \ref{prop:ConvexityWass}-$(c)$, $(d)$, it is clear that $\cG^{\blambda}_L(\bullet \,|\,\bnu,\rho)$ is an admissible surrogate for $\cJ^{\blambda}$. In the following theorem, we show iteratively minimizing the former generates sequences whose cluster points are critical for the median problem \eqref{eq:MMW2Med}. 

\begin{thm}[MM-scheme for approximating $W_2$-medians]
\label{thm:MM}
Suppose that $\bnu \in \cP_p^{\mathrm{ac}}(\R^d)^n$ for some $p \in (2,+\infty)$, fix $\blambda \in \Delta^n$, and let $\mu_0 \in \cP_p^{\mathrm{ac}}(\R^d)$. Then, the update rule 
\begin{equation}
\label{eq:MMUpdate}
\mu_{k+1} = \argmin_{\mu \in \cP_2(\R^d)} \cG^{\blambda}_L(\mu \,|\, \bnu,\mu_k) = \QrLWm_{2L}^{\blambda}(\bnu\,|\,\mu_k) 
\end{equation}
generates a bounded sequence $(\mu_k)_{k \in \N} \subset \cP_p^{\mathrm{ac}}(\R^d)$. Its cluster points $\mu^{\infty} \in \cP_p(\R^d)$ in the $W_2$-topology either satisfy $\mu^{\infty} = \nu_j$ for some $j \in \llbracket1,n\rrbracket$, or there exist optimal plans $\gamma_i^{\infty} \in \Gamma_o(\mu^{\infty},\nu_i)$ for which the following stationarity condition holds 
\begin{equation}
\label{eq:StationarityMM}
\sum_{i=1}^n \frac{\lambda_i}{W_2(\nu_i,\mu^{\infty})} \cdot \int_{\R^d} \langle x-y , \phi(x) \rangle \, \mathrm{d} \gamma_i^{\infty}(x,y) = 0
\end{equation}
for every $\phi \in C^0_b(\R^d,\R^d)$. In particular if $\mu^{\infty} \in\cP_2^{\mathrm{ac}}(\R^d)$, then 
%
%
the latter is critical for \eqref{eq:MMW2Med} in the sense that
\begin{equation}
\label{eq:StationarityMMBis}
\sum_{i=1}^n \lambda_i \cdot \nabla_{\mu} W_2(\nu_i,\mu^{\infty}) = 0
\end{equation}
in $L^2(\R^d;\R^d;\mu^{\infty})$.
\end{thm}

\begin{proof}
We showed in Theorem \ref{thm:QrLWm} that given any $\mu_k \in \cP_p^{\mathrm{ac}}(\R^d) \subset \cP_2(\R^d)$, the minimization problem in \eqref{eq:MMUpdate} admits a unique solution $\mu_{k+1} \in \cP_2(\R^d)$. Besides, it follows from the very definition of the iterates that $\cJ^{\blambda}(\mu_k) \leq \cJ^{\blambda}(\mu_0)$ for each $k \geq 1$ and 
\begin{equation*}
\begin{aligned}
\frac{1}{2L} W_2^2(\mu_k,\mu_{k+1}) & \leq \cG^{\blambda}_L(\mu_k \,|\, \bnu,\mu_k) - \cG^{\blambda}_L(\mu_{k+1} \,|\, \bnu,\mu_k) \\
& \leq \cJ^{\blambda}(\mu_k) - \cJ^{\blambda}(\mu_{k+1}).
\end{aligned}
\end{equation*}
By summing and telescoping, we then get that 
\begin{equation*}
\sum_{k=1}^{+\infty} W_2^2(\mu_k,\mu_{k+1}) \leq 2L \cdot \cJ^{\blambda}(\mu_0)
\end{equation*}
which entails in particular $W_2(\mu_k,\mu_{k+1}) \longrightarrow 0$. Now, since we assumed that $\bnu \in \cP_p^{\mathrm{ac}}(\R^d)^n$, it follows from the estimate \eqref{eq:MomentQrGm} that $\mu_{k+1} \in \cP_p^{\mathrm{ac}}(\R^d)$ with 
\begin{equation*}
\cM_p(\mu_{k+1}) \leq \max \Big\{ \cM_p(\mu_k) , \cM_p(\nu_1),\ldots, \cM_p(\nu_n)\Big\}. 
\end{equation*}
Whence, by induction, there exists a constant $C > 0$ which only depends on the magnitudes of the moments of $\mu_0 \in \cP_p^{\mathrm{ac}}(\R^d)$ and $\bnu \in \cP_2^{\mathrm{ac}}(\R^d)^n$ such that
\begin{equation}
\label{eq:MomentBoundMM}
\sup_{k \in \N} \cM_p(\mu_k) \leq C.
\end{equation}
Since we chose $p \in (2,+\infty)$, it follows e.g. from \cite[Proposition 7.1.5]{AGS} that the sequence $(\mu_k)_{k \in \N} \subset \cP_p^{\mathrm{ac}}(\R^d)$ is relatively compact in the $W_2$-metric, and thus admits converging subsequences.  

Let $\mu^{\infty} \in \cP_p(\R^d)$ be a cluster point of the MM-sequence attained along a subsequence $(\mu_{k_m})_{m \in \N}$. If $\mu^{\infty} \neq \nu_j$ for each $j \in \llbracket 1,n \rrbracket$, then there exists $\delta > 0$ and $M \in \N$ such that 
\begin{equation}
\min\bigg\{ \min_{i \in \llbracket 1,n \rrbracket} W_2(\nu_i,\mu_{k_m}) \, , \, \min_{i \in \llbracket 1,n \rrbracket} W_2(\nu_i,\mu_{k_m+1}) \bigg\}  \geq \delta 
\end{equation}
for every $m \geq M$, where we used the fact that both $(\mu_{k_m})_{m \in \N}$ and $(\mu_{k_m+1})_{m \in \N}$ converge to $\mu^{\infty}$. Recalling that the update rule \eqref{eq:MMUpdate} defines a quadratically regularized linearized Wasserstein median around the iterate $\mu_{k_m} \in \cP_p^{\mathrm{ac}}(\R^d)$, we infer from writing the first order optimality condition of the characterization \eqref{eq:QrGMlambda} thereof given in Theorem \ref{thm:QrLWm} that 
\begin{equation*}
\frac{1}{L} \Big( T_{\mu_{k_m}}^{\mu_{{k_m}+1}} - \Id \Big) + \sum_{i=1}^n \lambda_i \cdot \frac{T_{\mu_{k_m}}^{\mu_{{k_m}+1}} - T_{\mu_{k_m}}^{\nu_i}}{\big\| T_{\mu_{k_m}}^{\mu_{{k_m}+1}} - T_{\mu_{k_m}}^{\nu_i} \big\|_{L^2(\mu_{k_m})}} = 0
\end{equation*}
in $L^2(\R^d;\R^d;\mu_k)$ for each $m \geq M$. At this stage, we introduce for each $i \in \llbracket 1,n \rrbracket$ the transport plans $\bgamma_i^m := (\Id,T_{\mu_{k_m}}^{\nu_i},T_{\mu_{k_m}}^{\mu_{{k_m}+1}})_{\sharp} \mu_{k_m} \in \Gamma(\mu_k,\nu_i,\mu_{k+1})$, which by construction satisfy
\begin{equation*}
(\pi^1,\pi^2)_{\sharp} \bgamma_i^m \in \Gamma_o(\mu_{k_m},\nu_i) \qquad \text{and} \qquad (\pi^1,\pi^3)_{\sharp} \bgamma_i^m \in \Gamma_o(\mu_{k_m},\mu_{{k_m}+1}), 
\end{equation*}
and note then that 
\begin{equation}
\label{eq:StationarityMMProof}
\sum_{i=1}^n \int_{\R^{3d}} \bigg\langle \frac{1}{n L} \cdot (z-x) +  \lambda_i \cdot  \frac{z-y}{\big( \int_{\R^{3d}} \|z-y\|^2 \, \mathrm{d} \bgamma_i^m(x,y,z) \big)^{\nicefrac{1}{2}}} \, , \, \phi(x) \bigg\rangle \, \mathrm{d} \bgamma_i^m(x,y,z) = 0
\end{equation}
for every $\phi \in C^0_b(\R^d,\R^d)$. Since both $(\mu_{k_m})_{m \in \N}$ and $(\mu_{k_m+1})_{m \in \N}$ converge to the same limit $\mu^{\infty}$ in the $W_2$-topology, it follows from \cite[Proposition 7.1.3 and Proposition 7.1.5]{AGS} that there is a further subsequence (that we do not relabel) along which every $(\bgamma_i^m)_{m \in \N}$ converges in that same topology towards some $\bgamma_i^{\infty} \in \cP(\R^{3d})$ for which
\begin{equation*}
\gamma_i^{\infty} := (\pi^1,\pi^2)_{\sharp} \bgamma_i^{\infty} \in \Gamma_o(\mu^{\infty},\nu_i) 
\end{equation*}
and $(\pi^1,\pi^3)_{\sharp} \bgamma_i^{\infty} \in \Gamma_o(\mu^{\infty},\mu^{\infty}) = \{(\Id,\Id)_{\sharp} \mu^{\infty}\}$. Upon observing that the maps $(x,y,z) \in \R^{3d} \mapsto |\langle z-x , \phi(x) \rangle| \in \R$ and $(x,y,z) \in \R^{3d} \mapsto |\langle z-y , \phi(x) \rangle| \in \R$ are both uniformly integrable with respect to the sequences $(\bgamma^m_i)_{m \in \N}$, we may pass to the limit as $k \to +\infty$ in \eqref{eq:StationarityMMProof} (see e.g. \cite[Lemma 5.1.7]{AGS}) to obtain that 
\begin{equation*}
\begin{aligned}
& \sum_{i=1}^n \int_{\R^{3d}} \bigg\langle \frac{1}{n L} \cdot (z-x) +  \lambda_i \cdot \frac{z-y}{\big( \int_{\R^{3d}} \|z-y\|^2 \, \mathrm{d} \bgamma_i^{\infty}(x,y,z) \big)^{\nicefrac{1}{2}}} \, , \, \phi(x) \bigg\rangle \, \mathrm{d} \bgamma_i^{\infty}(x,y,z) \\
& = \sum_{i=1}^n \frac{\lambda_i}{\big( \int_{\R^{3d}} \|x-y\|^2 \, \mathrm{d} \gamma_i^{\infty}(x,y) \big)^{\nicefrac{1}{2}}} \cdot \int_{\R^d} \langle x-y, \phi(x) \rangle \, \mathrm{d} \gamma_i^{\infty}(x,y) \\
& = 0
\end{aligned}
\end{equation*}
because $x=z$ for $\bgamma^{\infty}_i$-almost every $(x,y,z) \in \R^{3d}$, which is precisely \eqref{eq:StationarityMM} since $(\pi^1,\pi^2)_{\sharp}\bgamma_i^{\infty} = \gamma_i^{\infty} \in \Gamma_o(\mu^{\infty},\nu_i)$ for each $i \in \llbracket 1,n \rrbracket$. 

In the case in which $\mu^{\infty} \in \cP_2^{\mathrm{ac}}(\R^d)$, it follows then from Theorem \ref{thm:Brenier} that $\gamma_i^{\infty} = (\Id,T_{\mu^{\infty}}^{\nu_i})_{\sharp} \mu^{\infty}$ for a uniquely determined Brenier map $T_{\mu^{\infty}}^{\nu_i} \in L^2(\R^d,\R^d;\mu^{\infty})$, wherefore 
\begin{equation*}
\sum_{i=1}^n \frac{\lambda_i}{W_2(\nu_i,\mu^{\infty})} \cdot (\Id - T_{\mu^{\infty}}^{\nu_i}) = 0 
\end{equation*}
by a routine density argument on continuous bounded maps. Recalling the expression \eqref{eq:GradWass} of the gradient of the squared Wasserstein distance whose existence is ensured by Theorem \ref{thm:WassGrad}, we end up with the stationarity condition \eqref{eq:StationarityMMBis}. 
\end{proof}

We conclude this section on MM approximations of Wasserstein medians by two elementary corollaries and remarks relating our approach to the literature.  

\begin{cor}[The case of bounded absolutely continuous inputs]
\label{cor:ClusterAC}
Suppose that $\nu_i = \rho_i \cdot \Leb$ for each $i \in \llbracket 1,n \rrbracket$ and let $\mu_0 := \rho_0 \cdot \Leb$ for some $(\rho_0,\ldots,\rho_n) \in L^{\infty}(\R^d,\R_+)^n$. Then, every cluster point of the MM-scheme is of the form $\mu^{\infty} = \rho^{\infty} \cdot \Leb$ with $\rho^{\infty} \in L^{\infty}(\R^d,\R_+)$.
\end{cor}

\begin{proof}
This comes as a direct consequence of the density bounds of Proposition \ref{prop:Density}-$(a)$ combined with the fact that those are stable under $W_2$-convergence.  
\end{proof}

\begin{cor}[Characterization of the cluster points as fixed-points]
\label{cor:FixedPoint}
Under the assumptions of Corollary \ref{cor:ClusterAC}, every cluster point $\mu^{\infty} \in \cP_2^{\mathrm{ac}}(\R^d)$ of the MM sequence is such that 
\begin{equation*}
\mu^{\infty} = \argmin_{\mu \in \cP_2(\R^d)} \cG^{\blambda}_L(\mu \,|\, \mu^{\infty}).    
\end{equation*}
In other words, every cluster point is a fixed-point of the iteration map. 
\end{cor}

\begin{proof}
From the moment estimate \eqref{eq:MomentBoundMM} and Corollary \ref{cor:ClusterAC}, we know that any sequence of iterates must lie inside a fixed compact set $\cK \subset \cP_2^{\mathrm{ac}}(\R^d)$, so that 
\begin{equation*}
\mu_{k+1} = \argmin_{\mu \in \cK} \cG^{\blambda}_L(\mu \,|\,\bnu,\mu_k)
\end{equation*}
for every $k \in \N$. Then, it follows from Berge's maximum theorem (see e.g. \cite[Theorem 17.31]{Aliprantis2006}) that the mapping 
\begin{equation*}
\rho \in \cK \mapsto \argmin_{\mu \in \cK} \cG^{\blambda}_L(\mu \,|\,\bnu,\rho)
\end{equation*}
is continuous provided that $(\mu,\rho) \in \cK^2 \mapsto \cG^{\blambda}_L(\mu \,|\,\bnu,\rho)$ is jointly continuous. Such a joint continuity follows from standard arguments akin to those of \cite[Theorem 3.1]{Alvarez2016}. Indeed, given converging sequences $(\mu_k)_{k \in \N} \subset \cP_2(\R^d)$ and $(\rho_k)_{k \in \N} \subset \cK$ with limits $\mu \in \cP_2(\R^d)$ and $\rho \in \cP_2^{\mathrm{ac}}(\R^d)$, the 3-plans defined by 
\begin{equation*}
\bgamma_k := (\Id,T_{\rho_k}^{\nu_i},T_{\rho_k}^{\mu_k})_{\sharp} \rho_k 
\end{equation*}
for each $i \in \llbracket 0,n \rrbracket$ are tight by construction (see e.g. \cite[Lemma 5.2.2]{AGS}), and by stability of optimal plans (see e.g. \cite[Proposition 7.1.3]{AGS}), its unique cluster point in the $W_2$-topology is 
\begin{equation*}
\bgamma = (\Id,T_{\rho}^{\nu_i},T_{\rho}^{\mu})_{\sharp} \rho. 
\end{equation*}
Combined with the standard convergence result e.g. of \cite[Lemma 5.1.7]{AGS}, this implies in particular that 
\begin{equation*}
W_2^2(\rho_k,\mu_k) ~\underset{k \to +\infty}{\longrightarrow} ~ W_2^2(\rho,\mu) \qquad \mathrm{and} \qquad W_{2,\rho_k}(\nu_i,\mu_k) ~\underset{k \to +\infty}{\longrightarrow} ~  W_{2,\rho}(\nu_i,\mu),  
\end{equation*}
and the continuity the iteration map follows. We conclude then by taking a converging subsequence and letting $k \to +\infty$ as in the proof of Theorem \ref{thm:MM}.
\end{proof}
    
\begin{rem}[Concerning the result of Theorem \ref{thm:MM}]
While the existence of cluster points result of Theorem \ref{thm:MM} may legitimately appear quite weak, it is worth stressing that they are somehow the best one can hope to get at this level of generality dealing with a nonconvex problem on a metric space, see \cite[Chapter 7]{Lange}, and is basically what is achieved in \cite{Alvarez2016,You2025}. A potential improvement would be to ensure that the sequence of iterates converge, which could be done e.g. by showing that the critical points of the $W_2$-median are isolated. Another more ambitious route would be to try and adapt the quantitative convergence analysis for MM-schemes developed in \cite{Bolte2016}. This seems however quite arduous at first glance, as the techniques therein extensively build on notions and results from real-algebraic geometry which do not seem to have any clear counterparts in Wasserstein spaces. Both endeavors are left for future works.   
\end{rem}

\begin{rem}[An MM perspective on \cite{Alvarez2016}]
It is possible to recast the fixed-point approach to Wasserstein barycenters developed in the above work as an MM-scheme, where
\begin{equation*}
\cJ^{\blambda}(\mu) := \sum_{i=1}^n \lambda_i \cdot W_2^2(\nu_i,\mu) \qquad \text{and} \qquad \cG^{\blambda}(\mu\,|\,\bnu,\rho) := \sum_{i=1}^n \lambda_i \cdot W_{2,\rho}^2(\nu_i,\mu). 
\end{equation*}
Indeed, as explained in the proof of Theorem \ref{thm:BreakdownWass2}, it holds that
\begin{equation*}
\Big( \sum_{i=1}^n \lambda_i \cdot T_{\mu_k}^{\nu_i} \Big)_{\sharp} \mu_k = \argmin_{\mu \in \cP_2(\R^d)} \cG^{\blambda}(\mu \,|\,\bnu,\mu_k). 
\end{equation*}
Although not formulated in this language, the qualitative results of \cite{Alvarez2016}, including proof of the well-posedness and continuity of the update map in Theorem 3.1, the decrease condition formulated in Proposition 3.3, and the convergence of subsequences towards fixed-point of the update in Theorem 3.6 all mirror the classical global MM convergence ansätze listed e.g. in \cite[Section 7.2]{Lange}. Under the working assumptions of their Theorem 3.6, it could be shown be similar arguments as those from Theorem \ref{thm:MM} above that every cluster point $\mu^{\infty} \in \cP_2^{\mathrm{ac}}(\R^d)$ for their scheme is critical for $\WBar_2^{\blambda}(\bnu)$, in the sense that 
\begin{equation*}
\sum_{i=1}^n \lambda_i \cdot \nabla_{\mu} W_2^2(\nu_i,\mu^{\infty}) = 0
\end{equation*}
in $L^2(\R^d,\R^d;\mu^{\infty})$. Reciprocally, it could be shown by adapting the arguments of Corollary \ref{cor:FixedPoint} above that cluster points $\mu^{\infty} \in \cP_p(\R^d)$ of an MM scheme would solve the fixed point equation
\begin{equation*}
\mu^{\infty} = \argmin_{\mu \in \cP_2(\R^d)} \cG^{\blambda}_L(\mu \,|\, \mu^{\infty}) = \Big( \sum_{i=1}^n \lambda_i \cdot T_{\mu^{\infty}}^{\nu_i} \Big)_{\sharp} \mu^{\infty}
\end{equation*}
which is exactly what \cite{Alvarez2016} achieves.
\end{rem}

\bibliographystyle{alpha}
{\footnotesize
\bibliography{references}}

\

\end{document}